\documentclass[11pt]{article}
\usepackage{amssymb}
\usepackage{latexsym}
\newcounter{MyCounter}
\newtheorem{theo}{Theorem}[section]
\newtheorem{pr}{Proposition}[section]

\newtheorem{ex}{Example}[section]
\newtheorem{re}{Remark}[section]
\newtheorem{co}{Corollary}[section]
\newtheorem{lm}{Lemma}[section]
\newcommand\proof{{\sl Proof. }}
\newcommand{\qed}{\hspace*{\fill}\hbox{$\Box$}\vspace{2ex}}
\newcommand{\Ad}{{\rm Ad}}
\newcommand{\ad}{{\rm ad}}

\newcommand{\be}{\begin{equation}}
\newcommand{\ee}{\end{equation}}
\newcommand{\bea}{\begin{eqnarray}}
\newcommand{\eea}{\end{eqnarray}}
\newcommand{\bean}{\begin{eqnarray*}}
\newcommand{\eean}{\end{eqnarray*}}
\newcommand{\benu}{\begin{enumerate}}
\newcommand{\eenu}{\end{enumerate}}
\newcommand{\edo}{\end{document}}

\newcommand{\barr}{\begin{array}}
\newcommand{\earr}{\end{array}}

\newcommand{\gep}{\varepsilon}
\newcommand{\eps}{\varepsilon}

\newcommand{\fg}{\mathfrak g}
\newcommand{\fd}{\mathfrak d}
\newcommand{\fh}{\mathfrak h}
\newcommand{\fa}{\mathfrak a}
\newcommand{\fb}{\mathfrak b}

\newcommand{\fp}{\mathfrak p}
\newcommand{\fz}{\mathfrak z}
\newcommand{\fsl}{{\mathfrak{sl}}}
\newcommand{\Dera}{\mathop{{\rm Der_a}}}
\newcommand{\fso}{{\mathfrak{so}}}

\newcommand{\fspin}{{\mathfrak{spin}}}
\newcommand{\fsu}{{\mathfrak{su}}}
\newcommand{\hol}{\mathfrak{hol}}

\newcommand{\Ker}{\mathop{{\rm Ker}}}
\newcommand{\Osc}{\mathop{{\rm Osc}}}
\newcommand{\osc}{\mathop{\mathfrak{osc}}}
\newcommand{\CC}{{\mathbb C}}
\newcommand{\TT}{{\mathbb T}}

\newcommand{\RR}{{\mathbb R}}

\newcommand{\NN}{{\mathbb N}}

\newcommand{\tH}{\tilde{H}}

\newcommand{\tX}{\tilde{X}}

\newcommand{\tal}{\tilde{\alpha}}

\newcommand{\tr}{{\rm tr}}
\newcommand{\Ric}{{\rm Ric}}
\newcommand{\Cliff}{{\rm Cliff}}

\newcommand{\SL}{{\rm SL}}

\newcommand{\SO}{{\rm SO}}

\newcommand{\Span}{\mathop{{\rm span}}}
\newcommand{\Spin}{{\rm Spin}}

\newcommand{\GL}{{\rm GL}}

\newcommand\ip{\mbox{$\langle\cdot \,,\cdot \rangle$}}
\begin{document}

\title{Doubly Extended Lie Groups - Curvature, Holonomy and Parallel Spinors}
\date{}
\author{Helga Baum and Ines Kath\footnote{The second author was supported by the
Max-Planck-Institut
 f\"ur Mathematik in Bonn, Germany.
}}
\maketitle
\begin{abstract}
In the present paper we 
study the geometry of doubly extended Lie groups with their natural biinvariant 
metric. We describe the curvature, the holonomy and the space of parallel spinors. 
This is completely done for all simply connected groups with biinvariant metric of 
Lorentzian signature $(1,n-1)$, of signature $(2,n-2)$ and of signature $(p,q)$, 
where $p+q\leq 6$. Furthermore, some special series with higher signature are 
discussed.
\end{abstract}
\tableofcontents
\section{Introduction}
The present paper is a step in a project which is aimed at a better understanding of 
the following  problems:
\begin{list}{\arabic{MyCounter}}
{\leftmargin0.4cm
\partopsep-0.4cm
\itemsep0cm}
\usecounter{MyCounter}
\item
 Which pseudo-Riemannian manifolds admit parallel spinors?
\item 
Describe (indefinite) metrics with indecomposable but non-irreducible holonomy 
representations.
\item
 Describe (pseudo-Riemannian) homogeneous spaces with non semi-simple isometry 
group. 
\end{list}
These problems are related to each other but are also of independent interest. 
The relation between parallel spinors and the holonomy representation is described 
by the well known holonomy principle, which states that the space of parallel 
spinors is isomorphic to the space of invariant elements of the holonomy 
representation on the spinor module. Homogeneous spaces with non semi-simple isometry 
group are an important class of examples of manifolds with indecomposable but 
non-irreducible holonomy representation.

Parallel spinor fields occur in several contexts. They are special solutions of the 
Dirac and the twistor equation of a semi-Riemannian spin manifold 
(\cite{Baum/Friedrich/ua:91}), they occur as a technical tool in constructing other 
kinds of special spinor fields, such as geometric Killing spinors (see 
\cite{Baer:93}, \cite{Baum2:89}, \cite{Kath:99}, \cite{Bohle:00}) and they are of 
interest in supergravity and string theories (e.g.~\cite{Figueroa:99}, 
\cite{Figueroa2:99}, \cite{Figueroa2:01}). Whereas in the pure Riemannian case the 
situation is quite well understood (\cite{Wang:89}, (\cite{Wang:95}, 
\cite{Moroianu/Semmelmann:00}) and the study is extended to other connections than 
the Levi-Civita connection (\cite{Moroianu:96}, \cite{Buchholz1:00}, 
\cite{Buchholz2:00}, \cite{Friedrich/Ivanov:01}) the question of the existence of 
parallel spinor fields even for the Levi-Civita connection in the Lorentzian and the 
general pseudo-Riemannian case is widely open and the occuring geometries are known 
only for irreducible manifolds (\cite{Baum/Kath:99}). The difficulty of the problem 
in the indefinite case is caused by the existence of pseudo-Riemannian manifolds 
with indecomposable but non-irreducible holonomy representations, which do not 
occur in the definite case. Contrary to the irreducible case these representations 
are not classified and the study of metrics realizing such holonomy representations 
is only in the beginning (see e.g. \cite{Berard-Bergery/Ikemakhen:93}, 
\cite{Figueroa:99}, \cite{Ikemakhen:96}, \cite{Ikemakhen1:99}, \cite{Boubel:00}). 
Already first examples show that parallel spinors are to expect mainly in the case 
of indecomposable but non-irreducible holonomy representations. This is seen already 
in the case of manifolds with transitive isometry group. Whereas there are no 
non-flat homogeneous Riemannian manifolds with parallel spinors, such spinors occur 
on all Lorentzian symmetric spaces with solvable transvection group 
(\cite{Baum1:00}). These Lorentzian symmetric spaces have non-irreducible abelian 
holonomy representation  (\cite{Cahen/Wallach:70}, \cite{Astrakhantsev:73}, 
\cite{Bieliavsky:92}). 

In several papers local normal forms of pseudo-Riemannian metrics with parallel 
spinors are constructed (\cite{Bryant:00}, \cite{Kath3:98}, \cite{Leistner:01}). 
We are mainly interested in {\it complete} examples. For that reason we are studying 
in the present paper the simplest case of complete pseudo-Riemannian spaces, the Lie 
groups with biinvariant pseudo-Riemannian metrics, calculate their holonomy, their 
curvature properties and the space of parallel spinors. Such groups are products of 
three types of Lie groups with biinvariant metrics: the semi-simple ones with 
Killing form, the abelian groups with flat metrics and the doubly extended Lie 
groups. Since parallel spinors occur only on scalar flat pseudo-Riemannian manifolds 
with 2-step nilpotent Ricci tensor, our attention was led to the third case, to that 
of doubly extended Lie groups. 

In particular we study the curvature, the holonomy and the existence of parallel 
spinors for the following classes of Lie groups:
\begin{enumerate}
\item
The cotangent bundle $T^*H$ of a simple Lie group $H$ with its natural biinvariant 
metrics of split signature (Theorem \ref{T4}). 
\item
Double extensions of abelian Lie groups by simple ones (Theorem \ref{T3}, Examples 
\ref{exhol}, \ref{su}, \ref{compact}).
\item  
Double extensions of abelian Lie groups by 1-dimensional ones (Theorem \ref{T5}),
\item
Special solvable Lie groups with maximal isotropic center (Theorem \ref{T6}).
\end{enumerate}
In particular, we classify all simply-connected Lie groups with biinvariant metric and describe its
parallel spinors up to dimension $n \leq 6$, for all dimensions $n$ in case of 
Lorentzian signature $(1,n-1)$ and in case of signature $(2,n-2)$ (Table 1 - Table 6).\\[-0.2cm]

%
\section{Lie groups with biinvariant metric}
\label{S2}
In this chapter we describe indecomposable Lie groups with biinvariant
indefinite metrics, their curvature and holonomy. Let $(D,d)$ be a Lie group
with biinvariant semi-Riemannian metric $d$ of signature $(p,q)$, let $\fd$ be 
the Lie algebra of $D$ and $B_{\fd}$ its Killing form. For the Levi-Civita
connection and the curvature tensor of $(D,d)$ one has the well-known formulas
\begin{eqnarray}
\nabla_X Y &=&{\textstyle \frac{1}{2}} [X,Y]\nonumber\\
R(X,Y)Z&=& - {\textstyle \frac{1}{4}} [[X,Y],Z]\\
\Ric |_{\fd \times \fd} &=& - {\textstyle \frac{1}{4}} B_{\fd}\label{E2}
\end{eqnarray}
for all left invariant vector fields $X,Y,Z \in \fd$. Let us denote by
$\ip_{\fd}$ the $\ad$-invariant scalar product of index
$(p,q)$ defined by the metric $d$ on the Lie algebra $\fd$. For the
Lie algebra $\hol (D,d)$ of the holonomy group of $(D,d)$  we have
\begin{eqnarray}
\hol (D,d) &=& {\rm span} \{R(X,Y) \mid X,Y \in \fd\} \nonumber \\
&=& \ad ([\fd , \fd])\subset \fso (\fd , \ip_{\fd} ) \,.\label{E3}
\end{eqnarray}
To shorten the notation we will call a Lie algebra with ad-invariant scalar product 
a {\em metric Lie algebra} and the corresponding Lie groups (equipped with the 
induced biinvariant metric) a {\em metric Lie group}. \\
The standard example for a metric Lie algebra is
a semi-simple Lie algebra $\fd$ with the metric $\ip_{\fd} =
B_{\fd}$ given by the Killing form. The corresponding Lie groups $(D,d)$
are semi-Riemannian Einstein spaces with negative scalar curvature $R= 
- \frac{1}{4} \cdot \dim \fd$. In particular, in this case there are no 
parallel spinors. The other extremal  case, an abelian Lie group
$D = \TT^r \times \RR^s$ with invariant semi-Riemannian metric is flat and 
its universal covering has the maximal possible number of linearely independent
parallel spinors. Let us now consider a third class of examples; the double
extensions of Lie groups. \\
Let $(\fg, \ip_{\fg})$ be a metric Lie algebra,
$(\fh , \ip_{\fh})$ a Lie algebra with ad-invariant symmetric
bilinear form (which can degenerate) and let $\pi : \fh \to \Dera(\fg , \ip_{\fg})$ 
be a Lie algebra homomorphism from $\fh$ into
 the Lie algebra of all antisymmetric derivations of $\fg$. Furthermore, we
denote by $\beta : \Lambda^2 \fg \to \fh^*$ the 2-cocycle
$$ \beta (X,Y)(H):= \langle \pi (H) X,Y \rangle_{\fg} \quad , \quad
X,Y \in \fg , H \in \fh $$

and by $\ad^*_{\fh} (H) : \fh^* \to \fh^*$ the coadjoint representation
of $\fh$.

Now we define a Lie algebra and a metric
structure on the vector space $\fd := \fh^* \oplus \fg \oplus \fh$. We will denote 
the resulting metric Lie algebra by $\fd_{\pi} (\fg , \fh)$. \\
The commutator is given by
\begin{eqnarray}\label{E4}
\lefteqn{[(\alpha , X, H) , (\tal , \tX , \tH)]_{\fd_\pi} :=}\\
&&(\beta(X, \tX)+\ad^*_{\fh} (H) \tal - \ad^*_{\fh} (\tH) \alpha \,, \,
[X, \tX]_{\fg} + \pi (H) \tX - \pi (\tH)X\,,\,
[H, \tH ]_{\fh})\,.\nonumber
\end{eqnarray}

The metric is defined by
\begin{equation}\label{E5}
\langle (\alpha , X,H),(\tal , \tX , \tH)\rangle_{\fd_\pi}:= \langle
H, \tH \rangle_{\fh} + \langle X,\tX\rangle_{\fg} + \alpha (\tH) +
\tal (H).
\end{equation}
The bilinear form $\ip_{\fd_\pi}$ is non-degenerate, 
$\ad_{\fd}$-invariant and of signature
$$ \mbox{signature} \, (\ip_{\fd_\pi})= \mbox{signature} \, 
(\ip_{\fg})+ (\dim \fh , \dim \fh) \ . $$
The Lie algebra $\fd_{\pi} (\fg , \fh)$ is called double extension of $\fg$ with 
respect to
$\pi$. It can be viewed as extension of the semi-direct product of
$\fg\rtimes_{\pi} \fh$ by the abelian Lie algebra $\fh^*$. Let $D_\pi$ be the simply 
connected Lie group with Lie algebra $\fd_{\pi} (\fg , \fh)$.
We will say that a metric Lie algebra is indecomposable if it is not the orthogonal sum of two
non-trivial metric Lie algebras. A metric Lie group is called indecomposable if its Lie algebra is
indecomposable.

\begin{theo}[Medina/Revoy] \label{T1} Let $(\fd , \ip_{\fd})$ be an 
indecomposable metric Lie algebra. Then we are in one of the
following cases:
\begin{enumerate}
\item $\fd$ is simple
\item $\fd$ is 1-dimensional
\item $\fd$ is a double extension $\fd =\fd_{\pi} (\fg , \fh)$ as
defined above, where $\fh$ is 1-dimensional or simple.
\end{enumerate}
\end{theo}
\proof See \cite{Medina/Revoy2:85}, \cite{Favre:87}, \cite{Bordemann:97}.
\qed

Since we can restrict ourselves to indecomposable Lie 
algebras, we can exclude the cases where $\pi \equiv 0$ and $\fg \not\equiv
0$, where all derivations $\pi(\fh) \subset \Dera (\fg , \ip_{\fg})$ are inner ones 
and where $\fg$ has a factor $\frak p$ with
$H^1 (\mathfrak p, \RR)=H^2(\mathfrak p,\RR)=0$ due to the following proposition:


\begin{pr}\label{P1} {\rm \cite{Figueroa:95}}
\begin{enumerate}
\item If $\pi \equiv 0$, then the Lie algebra $\mathfrak d_{\pi}(\fg,\fh)$ 
decomposes into an orthogonal 
direct sum of the Lie algebras $\fg$ and $\mathfrak h^* \rtimes_{\ad^*} \fh$.
\item If $\pi : \fh \to \Dera (\fg , \ip_{\fg})$
is given by inner derivations, i.e. if there is a homomorphism $\varphi:
\fh \to \fg$ such that 
$$ \pi (H) = \ad{}_\fg\varphi (H) \ , $$
then there exists an isometric Lie algebra isomorphism
$$\Phi : \fd_{\pi} (\fg , \fh) \longrightarrow \fg \times (\fh^* 
\rtimes_{\ad^*_{\fh}} \fh)$$
where the metric on the image is given by the product of the metrics $\ip_{\fg}$ and
$$\langle(\alpha, H),(\tal, \tH)\rangle_{\fh^* \rtimes_{\ad^*_\fh} \fh} = \langle
H, \tH \rangle_{\fh} + \langle \varphi (H) , \varphi (\tH)\rangle_{\fg}
+ \alpha (\tH) + \tal (H) $$
\item If $\fg$ has a factor $\mathfrak p$ with $H^1 (\fp ,\RR)= H^2 (\fp, \RR)=0$,
then its double extension $\fd_{\pi} (\fg , \fh)$ is decomposable.
\end{enumerate}
\end{pr}
\begin{pr} \label{P2} The Ricci curvature $\Ric^{D_{\pi}}$ and the
scalar curvature $R^{D_{\pi}}$ of the doubly extended Lie group $D_{\pi}$ satisfy
\begin{itemize}
\item[a)] $\Ric^{D_{\pi}} |_{\fg \times \fg} = \Ric^G |_{\fg \times \fg} =
- \frac{1}{4} B_{\fg}$\\[0.5em]
$\Ric^{D_{\pi}} |_{\fh^* \times \fd_{\pi}} =0$\\[0.5em]
$\Ric^{D_{\pi}} |_{\fh \times \fh} = \tr_{\fg} (\pi (\cdot ) \circ \pi (
\cdot )) + 2B_{\fh}$\\[0.5em]
$\Ric^{D_{\pi}} (X,H)=\tr_{\fg} (\ad_g (X) \circ \pi (H)) \quad \forall X \in \fg
, H \in \fh$
\item[b)] $R^{D_{\pi}} =R^G = - \tr_{\fg} B_{\fg}$.
\end{itemize}

In particular, we have
\begin{enumerate}
\item The Ricci tensor of $D_{\pi}$ is 2-step nilpotent if and only if the 
Ricci tensor of $G$  is 2-step nilpotent and 
\begin{eqnarray*}& \sum\limits^{n}_{j=1} \varepsilon_j B_{\fg} (X,X_j) \circ 
\tr_{\fg}
(\ad_{\fg} (X_j) \circ \pi (H))=0&\\
&\sum\limits_{j=1}^n \varepsilon_j \tr_{\fg} (\ad_{\fg} (X_j) \circ \pi
(H)) \cdot \tr_{\fg} (\ad_{\fg} (X_j) \circ \pi (\tH))=0&
\end{eqnarray*}
for all $H, \tH \in \fh$ and $X \in \fg$, where $X_1, \ldots , X_n$ is 
an ON-basis of $(\fg , \ip_\fg)$.
\item The Lie group $D_{\pi}$ is Ricci-flat if and only if\\[0.5em]
$B_{\fg} =0$\\[0.5em]
$\tr_{\fg} (\ad_g (X) \circ \pi (H))=0$\\[0.5em]
$\tr_{\fg} (\pi (H) \circ \pi (\tH)) + 2 B_{\fh} (H, \tH)=0$\\[0.5em]
for all $X \in \fg$ and $H, \tH \in \fh$.
\end{enumerate}
\end{pr}

\proof A direct calculation using formulas (\ref{E2}), (\ref{E4}) and (\ref{E5}).
\qed

Using formulas (\ref{E3}) and (\ref{E4}) we obtain the holonomy algebra of the 
doubly
extended Lie group $D_{\pi}$ (with its biinvariant pseudo-Riemannian metric).
\begin{pr}\label{P4} 
The holonomy algebra $\hol (D_{\pi})$ of the 
doubly extended Lie group $D_{\pi}$ is spanned by the following elements of $\mathfrak{so} (\fd_{\pi} (\fg , \fh))$
\begin{eqnarray*}
&&\left( 
\begin{array}{ccc}
0 & \beta (Y, \cdot) & -\ad_{\fh} (\cdot)^*\beta (Y_1 , Y_2)\\
0 & \ad_{\fg} Y & - \pi (\cdot )Y\\
0 & 0 & 0 \end{array} \right)   
,\quad
\left( 
\begin{array}{ccc}
\ad_{\fh} (H)^* & 0 & 0\\
0 & \pi (H) & 0\\
0 & 0 & \ad_{\fh} (H) \end{array} \right)  
,\\&&
\left( 
\begin{array}{ccc}
0  & 0 & \ad_{\fh}(\cdot)^* \alpha\\
0 & 0 & 0\\
0 & 0 & 0 \end{array} \right)    
,\quad
\left(
\begin{array}{ccc}
0 & \beta (X, \cdot ) & 0\\
0 & \ad_{\fg} (X) & - \pi (\cdot)X\\
0 & 0 & 0 \end{array}
\right),  
\end{eqnarray*}
where $Y = [Y_1 , Y_2 ]_{\fg}\in [\fg,\fg]_{\fg},\, H \in [\fh
, \fh ]_{\fh},\ \alpha \in \ad_{\fh}^* (\fh) (\fh^*), X \in \pi (\fh) \fg$
(all matrices are taken with respect
to the decomposition $\fd_{\pi} = \fh^* \oplus \fg \oplus \fh$).
\end{pr}


\section{Spin holonomy of double extensions}
Let $(V, \Phi)$ be an $n$-dimensional real vector space with scalar
product $\Phi$ of signature $(r,s)$. We denote by $\Cliff (V, \Phi)$ the
Clifford algebra generated by $(V, \Phi)$ with the relations
$$ x \cdot y + y \cdot x = -2 \Phi (x,y) {\bf 1} \quad , \quad x,y \in V \ . $$
If $n=2m$, then the complexified Clifford algebra $\Cliff^{\mathbb C} (V, \Phi)$ is isomorphic to
the algebra $\CC (2^m)$
of all complex $2^m \times 2^m$-matrices. In case $n=2m +1$ $\Cliff^{\mathbb C} (V, 
\Phi)$ is isomorphic to $\CC (2^m) \oplus \CC (2^m)$. For concrete
calculations we will use the following realization of these isomorphisms. Let
$e_1 , \ldots , e_n$ be an orthonormal basis of $(V, \Phi)$, let denote
by $\kappa_j$ the sign $\kappa_j = \Phi (e_j , e_j)= \pm 1$ and by $\tau_j$ the 
number
$$ \tau_j = \left\{ \begin{array}{ccl}
i & \mbox{if} & \kappa_j = \Phi (e_j , e_j)=-1\\[1em]
1 & \mbox{if} & \kappa_j =\Phi (e_j , e_j)=1
                    \end{array} \right. \ . $$
Furthermore, let
$$ U= \left( \begin{array}{cc}
i&0\\0&-i \end{array} \right) \quad , \quad
V= \left( \begin{array}{cc}
0&i\\i&0 \end{array} \right) \quad , \quad 
E= \left( \begin{array}{cc}
1&0\\0&1 \end{array} \right) \quad , \quad
T= \left( \begin{array}{cc}
0&-i\\i&0 \end{array} \right) \ . $$

Then in case $n=2m$ an isomorphism
$$ \Phi_{2m} : \Cliff^{\mathbb C} (V, \Phi) \longrightarrow \CC (2^m) $$
is given by\\
\parbox{14cm}{\begin{eqnarray}
\Phi_{2m} (e_{2j-1})&=& \tau_{2j-1} E \otimes \ldots \otimes E \otimes U
\otimes T \otimes \ldots \otimes T\nonumber\\[1em]
\label{E6}
\Phi_{2m} (e_{2j})&=& \tau_{2j} E \otimes \ldots \otimes E \otimes V
\otimes \underbrace{T \otimes \ldots \otimes T}_{j-1}
              \end{eqnarray}}
\hfill 

In case $n=2m+1$ we use the isomorphism
$$ \Phi_{2m+1} : \Cliff^{\mathbb C} (V, \Phi)  \to \CC (2^m) \oplus \CC (2^m) $$
given by\\
\parbox{14cm}{\begin{eqnarray}
\Phi_{2m+1} (e_{2j})&=& (\Phi_{2m} (e_j), \Phi_{2m} (e_j)) \quad j=1, 
\ldots , 2m\nonumber\\
\label{E7}
\Phi_{2m+1} (e_{n})&=& \tau_{n} (i T \otimes \ldots \otimes T, -i T \otimes 
\ldots \otimes T) \end{eqnarray}}
\hfill 

Let $\Spin (V, \Phi) \subset \Cliff (V, \Phi)$ be the Spin group of
$(V, \Phi)$. Let as above $m=[n/2]$. We denote by $\Delta_{r,s}$ the spinor 
representation of
$\Spin (V, \Phi)$ defined by
$$ \rho_{r,s} := \hat{\Phi}_{r+s} |_{\Spin (V, \Phi)} : \Spin
(V, \Phi) \to \GL (2^m,\CC) $$
where $\hat{\Phi}_{2m} := \Phi_{2m}$ and $\hat{\Phi}_{2m+1} = pr_1 \circ 
\Phi_{2m+1}$. Let us denote by $u(\varepsilon) \in \CC^2$ the vector 
$u(\varepsilon)= \frac{1}{\sqrt{2}} \left( {\textstyle \begin{array}{c} 1\\- 
\varepsilon i\end{array}} \right)$ , 
$\varepsilon = \pm 1$, and by 
$$ u(\gep_m , \ldots , \gep_1)=u(\gep_m) \otimes \dots\otimes u(\gep_1) \quad , 
\quad \gep_j = \pm 1 $$
the orthonormal basis of $\Delta_{r,s} = \CC^{2^m}$ with 
respect to the standard scalar product of $\CC^{2^m}$. The
Lie algebra $\fspin (V, \Phi)$ of ${\rm Spin} (V, \Phi)$ is the subalgebra
$$ \fspin (V, \Phi) = {\rm span} \{x \cdot y \mid x,y \in V\} \subset
\Cliff (V, \Phi) \ . $$
Let us denote by 
$$ \lambda : \Spin (V, \Phi) \to \SO (V, \Phi) $$
the 2-fold covering of the special orthogonal group defining the Spin 
group. Then the differential
$$ \lambda_* : \fspin (V, \Phi) \to \fso (V, \Phi) $$
satisfies
\begin{equation}\label{E8}
\lambda^{-1}_* (A) = \frac{1}{4} \sum\limits^n_{i=1} \kappa_i e_i \cdot A(e_i) 
\ , \quad A \in \fso (V, \Phi) \ . 
\end{equation}
Now, let $\fd := \fd_{\pi} (\fg, {\fh})$ be a doubly extended Lie
algebra as defined in Section \ref{S2}, where $(\fg, \ip_{\fg})$ is
$n$-dimensional of signature $(p,q)$ and $\dim \fh =r$. We choose a basis
$$ ( \underbrace{\alpha_r , \ldots , \alpha_1}_{\in \fh^*} , 
\underbrace{X_1 , \ldots , X_n}_{\in \fg} ,  
\underbrace{H_1 , \ldots , H_r}_{\in \fh} ) $$
such that 
$$ \langle X_i , X_j \rangle_{\fg} = \kappa_j \delta_{ij} \quad \quad \quad 
\kappa_j = \left\{
\begin{array}{cl} -1 & 1 \le j \le p\\
1& p<j \le p+q \end{array} \right.$$
$$ \langle H_i , H_j \rangle_{\fh} = c_j \delta_{ij} \quad , \quad
c_j \in \{ \pm 1,0 \}$$
$$ \alpha_i (H_j) = \delta_{ij} \ . $$
Using the orthonormal basis
\[
 e_{2i-1} = \frac{1}{\sqrt{2}} \Big(H_i - ( \frac{c_i}{2} +1) \alpha_i \Big),\qquad e_{2i}
=\frac{1}{\sqrt{2}} \Big(H_i - (\frac{c_i}{2} - 1) \alpha_i\Big),
\qquad i= 1, \ldots , r
\]
of $\fh^*\oplus\fh$ and 
$X_1 , \ldots , X_n$ of $\fg$
we obtain from (\ref{E8})
\begin{equation}\label{E9}
\lambda_*^{-1} (A)= \frac 14\sum\limits^r_{i=1} \Big(H_i A(\alpha_i) + \alpha_i 
A(H_i)-c_i\alpha_iA(\alpha_i)\Big)+ \frac14\sum\limits^n_{j=1} \kappa_j X_j \cdot 
A(X_j)
\end{equation}
for all $A \in \fso (\fd)$.
Furthermore, from the formulas (\ref{E6}) and (\ref{E7}) we obtain for the Clifford
multiplication in the spinor modul $\Delta_{p+r,q+r} = \Delta_{p,q} 
\otimes \Delta_{r,r}$
\begin{eqnarray}
\alpha_j \cdot (u \otimes u (\gep_r , \ldots , \gep_1)) &=&
\frac{1}{\sqrt{2}} (-1)^{j-1} \gep_1 \cdot \ldots \cdot\gep_{j-1} (\gep_j +1)
\nonumber \\
&& \qquad \qquad \cdot u \otimes u(\gep_r , \ldots , - \gep_j , \ldots \gep_1) 
\nonumber
\\[1em]\label{E10}
H_j \cdot (u \otimes u(\gep_r , \ldots , \gep_1)) &=& \frac{1}{\sqrt{2}}
(-1)^{j-1} \gep_1 \cdot \ldots \cdot\gep_{j-1} ((\gep_j -1)+
\frac{c_j}{2} (\gep_j +1)) \quad \quad \quad \mbox{} \\
&& \qquad \qquad\cdot u \otimes u( \gep_r , \ldots , - \gep_j, \ldots \gep_1) 
\nonumber \\[1em]
X \cdot (u \otimes u(\gep_r , \ldots , \gep_1))&=& (-1)^r \gep_1 \cdot
\ldots \cdot \gep_r (X \cdot u) \otimes u( \gep_r , \ldots , \gep_1) \nonumber
\end{eqnarray}
where $X \in \fg, u \in \Delta_{p,q} , j=1 , \ldots , r$.\\[0.2cm]
Let us denote by 
$$ \widetilde{\hol} (D_{\pi}):= \lambda_*^{-1} (\hol (D_{\pi})) \in 
\fspin (\fd)  $$
the holonomy algebra of the group $D_{\pi}$ represented in the Clifford algebra. 
Furthermore, we denote by
${\cal P}_{D_{\pi}}$ the space of parallel spinors of the group $D_{\pi}$. According 
to the holonomy principle this space is isomorphic to the elements of the spinor 
modul annihilated by the action of $\widetilde{\hol} (D_{\pi})$:
\[ {\cal P}_{D_{\pi}} \simeq \{ u \in \Delta_{p+r,q+r} \mid  \widetilde{\hol} 
(D_{\pi}) \cdot u = 0 \}\,.
\]

Combining formula (\ref{E9}) with Proposition \ref{P4} a direct calculation yields 
the
generators of the holonomy algebra
$ \widetilde{\hol} (D_{\pi})$.  


\begin{pr} \label{P5}
The holonomy algebra $\widetilde{\hol} (D_{\pi})$
is spanned by the elements 
$$ \sum\limits^r_{i,j=1} \langle \pi ([ H_i , H_j]_{\fh})Y_1 , Y_2 
\rangle_{\fg} \alpha_i \cdot \alpha_j 
+ 2 \sum\limits^r_{j=1} \pi (H_j) Y \cdot \alpha_j + 
\sum\limits^n_{j=1} \kappa_j X_j \cdot [Y, X_j]_{\fg} , $$
$$2 \sum\limits^r_{j=1} [H, H_j]_{\fh} \cdot \alpha_j + 2
\sum\limits^r_{j=1} \alpha_j([H, H_j]_{\fh}) \cdot {\bf 
1}+\sum\limits^r_{i,j=1}c_i\alpha_i([H,H_j])\alpha_i\cdot \alpha_j
- \sum\limits^n_{j=1} \kappa_j X_j \cdot \pi (H) X_j \  , $$
$$ \sum\limits^r_{i,j=1} \alpha ( [H_{i} , H_j ]_{\fh} ) \alpha_i
\cdot \alpha_j \quad , \quad 
 2 \sum\limits^r_{j=1} \pi (H_j) X \cdot \alpha_j + 
\sum\limits^n_{j=1} \kappa_j X_j \cdot [X, X_j ]_{\fg} \ , $$
where $Y = [Y_1 , Y_2]_{\fg} \in [\fg , \fg ]_{\fg} , \ H \in 
[\fh , \fh]_{\fh},\  \alpha \in {\rm ad}_{\fh} (\fh)^* \fh^* , \  X \in \pi
(\fh) \fg$.
\qed
\end{pr}

If the Lie algebra $\fh$ is abelian or semi-simple, then
$$ \sum\limits^r_{j=1} \alpha_j ([H, H_j]_{\fh})= {\tr}_{\fh}
({\rm ad} (H))=0 $$
for all $H \in \fh$. Hence, we obtain the following


\begin{theo}\label{T2}
Let $D_{\pi} = D_{\pi} (\fg , \fh)$ be a simply
connected metric Lie group, where $\fh$ is abelian or semi-simple, and let $G$ be the
simply connected metric Lie group corresponding to
$(\fg , \ip_{\fg})$. Let us denote by ${\mathfrak b}$ the
subspace
$$ {\mathfrak b}:= {\ad}_{\fg} (\pi(\fh) \fg)+\pi([\fh , \fh]_{\fh}) \subset \fso 
(\fg) . $$
Then the space ${\cal P}_{D_{\pi}}$ of parallel spinor fields on the double 
extension $D_{\pi}$
satisfies
$$ \dim {\cal P}_{D_{\pi}} \ge \dim (V_{\mathfrak b} \cap  V_{\hol (G)})
\ , $$
where for $\fa \subset \fso (\fg)$ the vector space $V_{\fa}$ is defined by
$$ V_{\fa} := \{ v \in \Delta_{p,q} \mid  \lambda^{-1}_* (\fa)\cdot v =0 \}\,.$$

In particular, if the Lie group $G$ admits a parallel spinor field 
defined by an element $v \in V_{\mathfrak{hol} (G)}$ such that $\lambda^{-1}_* 
({\mathfrak b})
\cdot v =0$, then $D_{\pi}$ admits a parallel spinor field.
\end{theo}

\proof Since $\hol (G)= {\rm ad}_{\fg} ([\fg , \fg]_{\fg})$ and 
$\fh$ is simple or 1-dimensional we obtain from Proposition \ref{P5} and the
formulas $(10)$ that
$$ \widetilde{\hol} (D_{\pi}) \cdot (v \otimes u(-1 , \ldots , -1))=0$$
for each $v \in V_{\hol (G)} \cap V_{\mathfrak b}$. 
The spinor $v \otimes u(-1 , \ldots , -1) \in \Delta_{p+r, q+r}$ defines
a parallel spinor field on $D_{\pi}$. 
\qed


\section{Double extensions by simple Lie groups}
Let $\fh$ be a simple Lie algebra. Then $[\fh, \fh]= \fh$ and ${\rm ad} 
(\fh)^* = \fh^*$. Therefore we have
$$[\fd , \fd]_{\fd} = \fh^* \oplus ([\fg , \fg]_{\fg} + \pi (\fh) \fg)
\oplus \fh \subset \fd \ . $$
Consequently, the holonomy algebra of $D_{\pi}$ is spanned by
$$ \left( \begin{array}{ccc}
0& \beta (X,\cdot\,)&0\\
0& \ad_{\fg} X & - \pi (\cdot)X\\
0&0&0 \end{array} \right) \ , \quad 
\left( \begin{array}{ccc}
\ad_{\fh} (H)^* &0&0\\
0& \pi (H) &0\\
0& 0 & \ad_{\fh}H \end{array} \right) \  , \quad
\left( \begin{array}{ccc}
0&0& \ad_{\fh} (\cdot)^* \alpha\\
0&0&0\\
0&0&0 \end{array} \right) $$
where $H \in \fh,\ \alpha \in \fh^*$ and $X \in [\fg , \fg]_{\fg} + \pi (\fh) 
\fg$.

Let us first consider the special case, that $(\fg , \ip_\fg)$
is an abelian Lie algebra with metric of signature $(p,q)$. Since $\fh$ is
simple the homomorphism
$ \pi : \fh \longrightarrow \fso (p,q) $
has no kernel. Hence $\fh$ can be viewed as simple subalgebra of $\fso (p,q)$.
>From Theorem \ref{T2} and Proposition \ref{P2} we obtain 

\begin{theo}\label{T3} Let $\fg$ be an abelian Lie algebra with scalar
product $\ip_{\fg}$ and let $\fh \subset \fso(\fg,\ip_\fg)$ be a simple subalgebra. 
We denote
by $\fd (\fg, \fh)$ the double extension given by $\pi = {\rm id}$ and
by $D(\fg, \fh)$ the simply connected Lie group corresponding to $\fd
(\fg , \fh)$. Then we have
$$ \dim {\cal P}_{D(\fg, \fh)} \ge \dim V_{\fh} \ . $$
The Ricci tensor of $D(\fg , \fh)$ is given
by
$$ \Ric (H, \tH)= {\rm Tr} (H \circ \tH)+ 2B_{\fh} (H, \tH) 
\not\equiv 0 \ . $$
All other components vanish.
\end{theo}

There is a whole string of examples where $\dim V_\fh\ge 1$. Let us mention the following ones.
\begin{ex}\label{exhol}
{\rm
In the paper \cite{Baum/Kath:99} we calculated the space $V_{\fh}$ for the simple 
Lie
algebras $\fh$ occurring in Berger's list of holonomy groups of irreducible 
pseudo-Riemannian manifolds and of $\fso^* (2p)$. It turns out that $\dim V_\fh\ge 1$
for $\fsu (p,q)\subset\fso(2p,2q)$, $\mathfrak{sp}(p,q)\subset\fso(4p,4q)$, $\fso^*
(2p)\subset\fso(2p,2q)$, $\fg_2\subset\fso(7)$, $\fg^*_{2(2)}\subset\fso(4,3)$, $
\fg^{\CC}_2\subset\fso(7,7)$, $ \fspin (7)\subset\fso(8)$, $\fspin (4,3)\subset\fso(4,4)$, and
$\fspin (7)^{\CC}\subset\fso(8,8)$.
}
\end{ex}
\begin{ex}\label{su}
{\rm
We consider the real irreducible representations of $\mathfrak{su}(2)$.
\begin{enumerate}
\item
\fbox{$\rho_k:\ \mathfrak{su}(2)\hookrightarrow \fso(2k+1),\ k\ge1$}\\[0.5ex]
The complexification $\rho_k^{\mathbb C}:\ \mathfrak{su}(2)
^{\mathbb 
C}=\mathfrak{sl}(2,\CC)\hookrightarrow \fso(2k+1,\CC)$ equals Sym$^{2k}\CC^2$, where 
$\CC^2$ is the standard representation of $\mathfrak{sl}(2,\CC)$. Consider
$H=\left(
\begin{array}{cc}
i&0\\0&-i
\end{array} 
\right)\,.$ Since the highest weight of $\rho_k^{\mathbb C}$ is $2k$ there exists an 
orthonormal frame $e_1,\dots,e_{2k+1}$ of $\RR^{2k+1}$ such that
\[\lambda_*^{-1}\rho_k(H)=e_1e_2+2e_3e_4+\dots+ke_{2k-1}e_{2k}\,.\]
The elements $u(\eps_k,\dots,\eps_1)$, $\eps_i=\pm 1$ for $i=1,\dots,k$ constitute a 
basis of weight vectors of the representation $(\lambda_*^{-1}\circ\rho_k)^{\mathbb 
C}$ of $\mathfrak{su}(2)^{\mathbb C}$ on $\Delta_{2k+1}$. More exactly we have
\[\lambda_*^{-1}\rho_k(H)u(\eps_k,\dots,\eps_1)=(i\eps_1+2i\eps_2+\dots+ki\eps_k)u(
\eps_k,\dots,\eps_1)\,.\]Let
\begin{eqnarray*}
N_0&:=&\#\{(\eps_k,\dots,\eps_1)\mid \eps_i=\pm 1\mbox{ for } i=1,\dots,k ; \ 
\eps_1+2\eps_2+\dots+k\eps_k =0\}\\
N_2&:=&\#\{(\eps_k,\dots,\eps_1)\mid \eps_i=\pm 1\mbox{ for } i=1,\dots,k ; \ 
\eps_1+2\eps_2+\dots+k\eps_k =2\}\ .
\end{eqnarray*}
The decomposition of $(\lambda_*^{-1}\circ\rho_k)^{\mathbb C}$ into irreducible 
$\mathfrak{sl}(2,\CC)$-representations yields
\[\dim V_{\rho_k( \mathfrak{su}(2)^{\mathbb C})}=N_0-N_2\,.\]
Consequently, there exists a biinvariant metric of signature $(3,2k+4)$ on \\
 $D(\RR^{2k+1},\rho_k(\mathfrak{su}(2)))$ with at least $N_0-N_2$ parallel spinors. 
For example, we have
\begin{center}
\begin{tabular}{|c||c|c|c|c|c|c|c|c|c|c|}\hline
\mbox{}&&&&&&&&&&\\[-0.5em]
$k$ & $3$& $4$ &$7$&$8$ &$11$ & $12$ &$15$ &$16$ &  $19$ &$20$   \\[0.5em]
\hline
\mbox{}&&&&&&&&&&\\[-0.5em]
$N_0-N_2$&
  $1$& $0$&$0$&$1$&$1$&
$1$&
$ 3$&$1$&$5$&$12$\\[0.5em]
\hline
\end{tabular}
\end{center}
If $k\equiv 1,2$ mod $4$, then $N_0-N_2=0$. 

\item
\fbox{$\sigma_k:\ \mathfrak{su}(2)\hookrightarrow \fso(4k),\ k\ge1$}\\[0.5ex]
The representation $\sigma_k$ equals $($Sym$^{2k-1}\CC^2)_{\mathbb R}$, where 
$\CC^2$ is the standard representation of $\mathfrak{su}(2)$. Since the highest 
weight of Sym$^{2k-1}\CC^2$ equals $2k-1$ there exists an orthonormal frame 
$e_1,\dots,e_{4k}$ of $\RR^{4k}$ such that
\[\lambda_*^{-1}\sigma_k(H)=\frac 
12\Big((e_1e_2-e_3e_4)+3(e_5e_6-e_7e_8)+\dots+(2k-1)(e_{4k-3}e_{4k-2}-e_{4k-1}e_{4k}
)\Big)\,.\]
Hence, $u(\eps_k,\dots,\eps_1)$, $\eps_i=\pm 1$ for $i=1,\dots,k$ is a basis of 
weight vectors of $(\lambda_*^{-1}\circ\sigma_k)^{\mathbb C}$ of 
$(\Delta_{4k},\lambda_*^{-1}\sigma_k)^{\mathbb C}$. And with
\begin{eqnarray*}
N_0':=\#\{(\eps_{2k},\dots,\eps_1)&\mid& \eps_i=\pm 1\mbox{ for } i=1,\dots,k ; \\&&
 (\eps_1-\eps_2)+3(\eps_3-\eps_4)+\dots+(2k-1)(\eps_{2k-1}-\eps_{2k}) =0\}\\
N_2':=\#\{(\eps_{2k},\dots,\eps_1)&\mid& \eps_i=\pm 1\mbox{ for } i=1,\dots,k ; \\&&
 (\eps_1-\eps_2)+3(\eps_3-\eps_4)+\dots+(2k-1)(\eps_{2k-1}-\eps_{2k}) =4\}\
.
\end{eqnarray*}
we obtain
\[\dim V_{\sigma_k( \mathfrak{su}(2))}=N_0'-N_2'\,.\]
Furthermore, since $\sigma_k(\mathfrak{su}(2))\subset 
\mathfrak{su}(2k)\subset\fso(4k)$ we have $N_0'-N_2'\ge 2$ (see Example 
\ref{exhol} and \cite{Baum/Kath:99}).
Consequently, there exists a biinvariant metric of signature $(3,4k+3)$ on  
$D(\RR^{4k},\sigma_k(\mathfrak{su}(2)))$ with at least $N_0'-N_2'\ge 2$ parallel 
spinors. For example, we have

\begin{center}
\begin{tabular}{|c||c|c|c|c|c|c|c|c|c|c|}\hline
\mbox{}&&&&&&&&&&\\[-0.5em]
$k$ & $1$& $2$ &$3$ & $4$ &$5$ &$6$ &  $7$ &$8$ &$9$ &$10$ \\[0.5em]
\hline
\mbox{}&&&&&&&&&&\\[-0.5em]
$N_0'-N_2'$&
  $2$& $3$&$4$&
$5$&
$ 8$&$11$&$16$&$29$&$50$&$94$\\[0.5em]
\hline
\end{tabular}
\end{center} 
\end{enumerate}
}
\end{ex}
\begin{ex}\label{compact}
{\rm
Let $\fh$ be a compact simple Lie algebra and $(V,\rho)$ a representation of 
$\fh^{\mathbb C}$. Then there exists a hermitian product $\ip$ on $V$ such that
\[\rho(\fh)\subset \mathfrak{su}(V,\ip)\subset \fso(V_{\mathbb R}, {\rm Re}\ip)\,.\]
Consequently, 
\[\dim V_{\rho(\fh)}\ge 2\]
(see Example 
\ref{exhol} and \cite{Baum/Kath:99}). Hence, $D(V_{\mathbb R}, \rho(\fh))$ is a Lie group 
with biinvariant metric of signature $(\dim \fh,\dim V_{\mathbb R}+\dim\fh)$ and 
admits at least two parallel spinors.
}
\end{ex}

Finally, we consider the case, where $\fg =\{0\}$ and $\fh$ is simple. We equip 
$\fh$ with a multiple $c \cdot B_{\fh}$ of its Killing form, $c \in \RR$. Then 
the double extension $\fd$ is isomorphic to the semi-direct product $\fd =
\fh\;_{\ad^*}  \hspace{-0.3pt} \ltimes \fh^*$ of $\fh$ with the abelian Lie algebra 
$\fh^*$. Let
$D$  be the corresponding Lie group $D= H \;_{\Ad^*}\hspace{-0.3pt}\ltimes\fh^*$.
$D$ is isomorphic to the cotangent bundle $T^*H$ of the simple Lie group $H$. If 
$\fh$ is equipped with the metric $c\cdot B_{\fh}$, we denote the corresponding 
metric Lie group with $T^*H_c$. For different parameters $c_1$ and $c_2$ the groups 
$T^*H_{c_1}$ and $T^*H_{c_2}$ are not isometric isomorphic. Indeed, consider a Levi 
decomposition $\mathfrak{rad}\oplus{\mathfrak s}$ of $\fh\;_{\ad^*}  \hspace{-0.3pt} 
\ltimes \fh^*$. Then the (uniquely determined) radical $\mathfrak{rad}$ is equal to 
$\fh^*\subset \fh\;_{\ad^*}  \hspace{-0.3pt} \ltimes \fh^*$ and by the Theorem of 
Malcev-Harish-Chandra there exists an element $u\in\fh^*\subset D$ 
such that $\Ad(u)\fh={\mathfrak s}$. Since the bilinear form $\ip_\fd$ on 
$\fh\;_{\ad^*}  \hspace{-0.3pt} \ltimes \fh^*$ is Ad-invariant we can recover 
$cB_\fh$ (and therefore $c$) from $\ip_\fd$.


\begin{theo}\label{T4}
Let $H$ be a simply connected simple Lie group. Then 
there is an exactly one-dimensional space of parallel spinor fields on the 
cotangent bundle $T^* H_c$. The Ricci tensor of $T^*H_c$  is given by
$$ {\rm Ric}^{T^*H} (Z, \tilde{Z})=2B_{\fh} (Z, \tilde{Z}) \ , $$
where $Z, \tilde{Z} \in \fh$. All other components are zero.

\end{theo}
\proof Let us first prove that the only subspaces of $\fd$ which are invariant under 
the action of the holonomy group are $\{0\}$, $\fh^*$, and $\fd$. Since $\ad_\fd 
(\fh)$ is a subspace of the holonomy algebra which acts on $\fh^*\subset \fd$ by the 
co-adjoint representation and $\fh$ is simple any invariant subspace of $\fh^*$ is 
equal to $\{0\}$ or $\fh^*$. Let now $U\subset\fd$ be an arbitrary invariant 
subspace. Then also $U\cap\fh^*$ is invariant. Hence $U\cap \fh^*=0$ or 
$\fh^*\subset U$. In the latter case $U/\fh^*$ is an invariant subspace of 
$\fd/\fh^*$. Again, $\ad_\fd (\fh)$ acts by adjoint representation and has no other 
invariant subspaces than $\{0\}$ and $\fd /\fh^*$. Hence, in this case $U=\fh^*$ or 
$U=\fd$. Now consider the case $U\cap \fh^*=0$ and let $h+h^*$ be in $U$, $h\in H$, 
$h^*\in H^*$. Since $U$ is invariant,  
$\ad_\fd(\alpha)(h+h^*)\in U$ for all $\alpha\in\fh^*$. 
Hence, $\ad_\fh(h)^*\alpha\in \fh^*\cap U$. Because of $U\cap \fh^*=0$ we obtain 
$\ad_\fh(h)^*\alpha=0$ for all $\alpha\in\fh^*$. Hence, $\ad(h)=0$, and since $\fh$ 
is simple, $h=0$. Consequently, $h^*\in U\cap\fh^*=0$. This implies $U=0$.

Proposition \ref{P5} shows that $u(-1,\dots,-1)$ is in $V_{\mathfrak{hol}(D)}$ and, 
hence, it defines a parallel spinor field $\psi$ on the 
cotangent bundle $T^* H$. 
Now let $\psi'$ be a further parallel spinor on $T^*H$. Let it be given by  $v\in 
V_{\mathfrak{hol}(D)}$. We define
\[U=\{Z\in TT^*H\mid Z\cdot\psi'=0\}\,.\] Then $U$ is a totally isotropic parallel 
distribution on $T^*H$ and, hence, given by a totally isotropic subspace of $\fd$ 
which is invariant under the action of the holonomy algebra. Consequently, this 
subspace can only be $0$ or $\fh^*$. If it is equal to $\{0\}$ then 
$\Ric(X)\cdot\psi'=0$ implies $\Ric(X)=0$ which contradicts $\Ric^{D} |_{\fh \times 
\fh} = 2B_{\fh}\not=0$.  Hence it must be equal to $\fh^*$ and, thus,
\[\fh^*\cdot v =0\,.\] Now (\ref{E10}) implies that $v$ is a multiple of 
$u(-1,\dots,-1)$. Consequently, $\psi'$ is a multiple of $\psi$.
\qed

%
          \section{Double extensions by 1-dimensional Lie groups}
\label{S4}
%
In this chapter we consider the case, that $\fh$ is 1-dimensional. Let $H$ be a generator of $\fh$
and $\langle H,H\rangle_\fh=c$. Then
the double extensions of the Lie algebra $(\fg, \ip_\fg)$ are
described by antisymmetric derivations $A \in \Dera (\fg, \ip_{\fg})$:
$$\fd_A :=\fd_A(\fg,\RR):= \RR \alpha \oplus \fg \oplus \RR H $$
with the commutator
\begin{eqnarray}
\label{EC1}
[\alpha,\,\cdot\;]_{\fd_A} &=& 0\nonumber\\ 
\ [X,Y]_{\fd_A} &=& \langle AX, Y \rangle_{\fg} \alpha + [X,Y]_{\fg}\\  
\ [H,X]_{\fd_A} &=& AX\nonumber
\end{eqnarray}
for all $X,Y \in \fg$. The invariant metric on $\fd_A$ defined by (\ref{E5}) equals
$$\ip_c := \left( \begin{array}{ccc}
0&0&1\\
0& \ip_{\fg} & 0\\
1&0&c \end{array} \right) \quad , \quad c \in \RR \ . $$
It is not hard to prove that all $(\fd_A, \ip_c)$ for $c\in\RR$ are isomorphic as metric Lie
algebras. Hence, we will assume that $c=0$.
 The centre of $\fd_A$ 
is
\begin{equation} \label{Ec}\fz (\fd_A)= \left\{ \begin{array}{lll}
\RR \alpha \oplus (\ker A \cap \fz (\fg)) &\mbox{ if }& A\not\equiv 0 \mbox{ mod } 
\ad(\fg)\\ 
\RR \alpha \oplus \fz (\fg)\oplus \RR(H-X_0) &\mbox{ if }& A=\ad(X_0)\ .
\end{array}\right. \end{equation}
If we suppose that $\fd_A$ is indecomposable, then  $A  \not\equiv 0$ mod $\ad(\fg)$ 
and $\ker A \cap \fz (\fg)$ is totally
isotropic. The following proposition shows that it suffices to consider outer antisymmetric
derivations
$ A\in{\rm Out}_{\rm a} (\fg, \ip_{\fg})= \Dera (\fg , \ip_{\fg} )/ {\rm ad} (\fg) \ . 
$
\begin{pr}\label{P6} {\rm \cite{Favre:87}}
 Let $(\fg , \ip_{\fg} )$ be a metric Lie algebra
 and $A, \hat{A} \in \Dera (\fg , \ip_{\fg})$ two antisymmetric derivations on 
$\fg$. Then there exists
an isometric Lie algebra isomorphism
$$\Psi : \fd_A (\fg, \RR )\longrightarrow \fd_{\hat{A}} (\fg , \RR)$$
if and only if there exists $\lambda_0 \in \RR \backslash \{ 0 \}, \ T_0 \in \fg$
and $\varphi_0 \in {\rm Aut} (\fg , \ip_{\fg})$ 
such that
$$ \varphi_0^{-1} \hat{A} \varphi_0 = \lambda_0 A + \ad (T_0) \ . $$
\end{pr}

We denote by $D_A(\fg)$ a metric Lie group with 
Lie algebra
$\fd_A = \fd_A (\fg , \RR) \ . $

Now, let us consider the special case, that 
 $(\fg , \ip_{\fg} )$ is an abelian
pseudo-Euclidean Lie algebra of signature $(p,q)$ and $A \in \fso (p,q)$ is a
non-vanishing antisymmetric endomorphism on $\fg$. 
In this case we denote the Lie algebra $\fd_A=\fd_A(\fg,\RR)$ by the symbol 
$\underline{A}(p,q)$ and the corresponding  simply connected Lie group by $A(p,q)$. 

\begin{theo} \label{T5} Let $\underline A(p,q)$ be the double extension of an abelian metric Lie
algebra of signature $(p,q)$ by $\RR$ defined by $A\in\fso(p,q)$ and let $A(p,q)$ be the simply
connected Lie group with Lie algebra $\underline A (p,q)$. The group $A(p,q)$ is indecomposable if
and only if $\Ker A$ is totally isotropic. Furthermore, $A(p,q)$ is solvable, scalar flat and has a
2-step nilpotent Ricci tensor, where
\begin{eqnarray*}
A(p,q) \quad  \mbox{is Ricci-flat} & \Leftrightarrow & \tr_{\fg} 
(A^2)=0\\
A(p,q) \quad \mbox{is flat} & \Leftrightarrow & A^2 =0\,.
\end{eqnarray*}
The holonomy of $A(p,q)$ is abelian with Lie algebra
$$ \hol (A(p,q)) = \left\{ \left( \begin{array}{ccc}
0& (A^2 X)^* &0\\
0&0& -A^2 X\\
0&0&0
\end{array} \right)\ \Big|\ X\in \fg \right\}\subset \fso (p+1,q+1) \ . $$
If the space $A^2\fg$ is totally isotropic, then the space ${\cal 
P}_{A(p,q)}$ of parallel spinor fields on $A(p,q)$ has dimension  
$\;2^{[\frac{p+q}{2}]}(1+2^{-\rho })\;$ , where $\rho $ is the dimension of 
$A^2\fg$. 
In case $A^2\fg$ is not totally isotropic, the dimension of ${\cal P}_{A(p,q)}$ is
$\;2^{[\frac{p+q}{2}]}$.
\end{theo}
\proof
The claimed indecomposability is not hard to prove. The assertion about the Ricci tensor follows
from Proposition \ref{P2}.
The holonomy algebra $ \widetilde{\hol} (A(p,q))$ in this case is 
$$ \widetilde{\hol} (A(p,q)) = \{ A^2 X \cdot \alpha \quad | \quad X  \in
\fg \} \ . $$
Formulas $(\ref{E10})$ show that $ w \otimes u(1)+ v \otimes u(-1)$
belongs to the annihilator of $\widetilde{\hol} (A(p,q)) \subset \fspin (\fd_A)$, 
if and only if
$$ A^2 X \cdot w =0 \quad \quad \forall X \in \fg \ . $$
If ${\rm Im}
A^2$ is not totally isotropic, then $w=0$, otherwise $w$ lies in the annihilator of 
$A^2\fg$, which has dimension $\;2^{[\frac{p+q}{2}]-\rho}\;$, if $\rho$ is the 
dimension of $A^2\fg$. \qed
\begin{re}
{\rm 
The group $A(p,q)$ is nilpotent if and only if $A$ is nilpotent. In this case 
$A(p,q)$ is an 
Einstein space of scalar curvature 0 (non-flat if $A^2 \not= 0$). The
nilpotent Lie algebras $\underline{A}$ are classified, cf. \cite{Favre:87}.
}
\end{re}
%
  \section{Indecomposable simply connected Lie groups with biinvariant Lorentzian
metric}
%

In this section we recall the classification of simply connected Lie groups with 
biinvariant Lorentzian metric and we determine the space of parallel spinors on these groups. If
$D$ is a simply connected 
simple Lie group with 
biinvariant Lorentzian metric, then it is isomorphic to
$\widetilde{\SL} (2, \RR)$ since  $\mathfrak{sl}(2, \RR)$ is the
only simple Lie algebra for which the Killing form has Lorentzian signature. 

\begin{pr}{\rm \cite{Medina:85}} Each indecomposable non-simple metric Lie algebra of Lo\-rent\-zian signature
$(1,n-1)$ is isomorphic to exactly one of the double extensions
\[\underline A_{\underline \lambda}(0,2m),\ \underline{\lambda}=(\lambda_1,\dots,\lambda_m),\ 
\lambda_1=1\le\lambda_2\le\dots\le\lambda_m\,,\] 
where $n=2m-2$ and 
\begin{equation}\label{EAlambda}
A_{\underline \lambda}= \left( \begin{array}{ccc}
\Lambda_1 && 0\\
& \ddots &\\
0&& \Lambda_m
             \end{array}
\right) \ ,\ \mbox{  } 
\Lambda_r=\left(\begin{array}{cc} 0&-\lambda_r\\ \lambda_r &0 \end{array}\right)\,. 
\end{equation}
\end{pr}
\proof
By Theorem \ref{T1} each non-simple indecomposable Lie algebra with ad-invariant 
Lo\-rentzian scalar product   
is the double extension $\fd_A (\fg ,\RR)$ of a  Euclidean
Lie algebra $(\fg,\ip_\fg)$ by $\RR$ defined by an antisymmetric 
derivation $A \in \Dera (\fg)$. Since $\ip_\fg$ is positive definite, $\fg$ splits orthogonally
into an 
abelian ideal $\fa$ and a semi-simple ideal $\fp$. By Proposition \ref{P1} (3.), 
$\fg$ cannot have a semi-simple factor, hence $\fp=0$. Therefore, $\fg$ is abelian 
and $A \in \fso(\fg)$ is a bijective antisymmetric map. In particular, the dimension 
of $\fg$ is even. By Proposition \ref{P6} $(\fd_A (\fg, \RR),\ip_\fd)$ and 
$(\fd_{\hat A} (\fg, \RR),\ip_\fd)$ are isomorphic if and only if there is a map 
$\varphi_0\in O(\fg,\ip_\fg)$ and a real number $\lambda_0\not= 0$ such that 
$\varphi_0^{-1}\hat A\varphi_0=\lambda_0A$. If 
$\dim(\fg)=2m$ each such class of maps in $\fso(\fg)$ is represented by exactly one 
of the maps
$A_{\underline \lambda}$
for $\underline{\lambda}=(\lambda_1,\dots,\lambda_m),\ 
\lambda_1=1\le\lambda_2\le\dots\le\lambda_m$.
\qed

The simply connected Lie group $A_{\underline \lambda}(0,2m)$ is also known as oscillator group and
denoted by $\Osc(\underline{\lambda})$. By construction this group is a semi-direct product of
$\RR$ with the Heisenberg group (see \cite{Medina:85}).

The group $\widetilde{\SL} (2,\RR)$ with its Killing form is an Einstein space
of negative scalar curvature $R= - \frac{3}{4}$. The oscillator group
$\Osc(\underline{\lambda})$ has a 2-step nilpotent Ricci 
tensor which is given by
$$ \Ric (H,H) = -\sum\limits^m_{i=1} \lambda^2_i \ , $$
All other components of the Ricci tensor
are zero (see Proposition \ref{P2}).

It follows from Theorem 
\ref{T5} that the $2m+2$-dimensional oscillator groups $\Osc(
\underline{\lambda})$, have $2^m$-linearly independent parallel spinor
fields.

We summarize the results of this section in the following table. 

\newpage

{\bf Table 1: Isometry classes of indecomposable simply connected Lie groups with biinvariant
Lorentzian metric and number of parallel spinors}
\begin{center}
\begin{tabular}{|l|c|c|}
\hline &&\\[-0.5ex]
$G$ & $\dim G$ & $\dim {\cal P}_G$  \\[2ex] \hline \hline &&\\[-0.5ex]
$(\widetilde{\SL} (2, \RR), c B_{\fsl (2, \RR)}), c \in \RR^+$ & $3$ & $0$ \\[2ex] \hline &&
\\[-0.5ex]
$ Osc (\lambda_1, \ldots , \lambda_m),\, 0 < \lambda_1 = 1 \le \lambda_2 \le \ldots 
\le \lambda_m, \,
m \ge 1 $ & $2m+2$ & $2^m$ \\[2ex] \hline
\end{tabular}
\end{center}

%
   \section{Indecomposable simply connected Lie groups with biinvariant 
             metric of signature $(2,n-2)$}
%

Now we turn to metric Lie groups of index 2. First we will describe the structure of such groups.
In principle all these groups have to be included in the classification of symmetric spaces with
index 2 due to Cahen and Parker (\cite{Cahen/Parker:80}, see also \cite{Neukirchner:02}). However, it is
not an easy task to identify the indecomposable Lie groups among these spaces.  Moreover, there are
metric groups which are isomorphic as symmetric spaces but not as metric groups. Also the notion of
indecomposability is not the same for groups and for symmetric spaces. Therefore, we prefer to use
the Theorem of Medina and Revoy. Clearly, $\widetilde{SL}(2,\RR)$ is the only Lie group with
biinvariant metric of index 2 which is simple. The non-simple indecomposable ones can be
characterized as follows.
\begin{pr}{\rm (\cite{Kath/Olbrich:02})} \label{P7.1}
If $(\fd,\ip)$ is a non-simple indecomposable Lie algebra with ad-invariant metric of signature
$(2,n-2)$, then $(\fd,\ip)$ is isomorphic to a double extension $\fd_A(\fg,\RR)$, where
\begin{enumerate}
\item
$\fg$ is an abelian Lorentzian Lie algebra and $A$ is an antisymmetric
endomorphism on $\fg$ with totally isotropic kernel, or
\item $\fg$ is an oscillator algebra 
$$ \fg = \fd_{A_0} (\fg_0 , \RR)= \RR \alpha_1 \oplus \fg_0 \oplus \RR H_1 , $$
where $A_0 \in \fso (\fg_0)$ is a bijective antisymmetric map on the
abelian Euclidean Lie algebra $\fg_0$. With respect to the decomposition
$\RR \alpha_1 \oplus \fg_0 \oplus \RR H_1$ the antisymmetric derivation
$A \in \Dera (\fg)$ is given by
$$ A= \left( \begin{array}{ccc}
0&0&0\\
0&U_1&0\\
0&0&0
\end{array}
\right) , $$
where $U_1 \in \fso (\fg_0)$ commutes with $A_0$ and there is no
decomposition $\fg_0 = \fg_0' \oplus \fg_0''$ of $\fg_0$ into orthogonal
$A_0$-invariant subspaces such that
$$ U_1 = \left( \begin{array}{cc}
t' A'_0 &0\\
0& t'' A''_0    \end{array} \right)
:
 \fg_0' \oplus \fg_0'' 
\rightarrow
\fg_0'  \oplus \fg_0'' , $$
where $A'_0=A_0|_{\fg_0'},\ A''_0=A_0|_{\fg_0''}$ and $t', t'' \in \RR$, or \\
\item 
$\fg$ is the direct sum of a one-dimensional Lie algebra and an oscillator algebra $\fd_{A_0}(
\fg_0 , \RR )$ with a bijective antisymmetric map $A_0 \in \fso (\fg_0)$, i.e.
$$ \fg = \RR\oplus\fd_{A_0} (\fg_0 , \RR)= \RR \oplus \RR \alpha_1
\oplus \fg_0 \oplus \RR H_1 . $$
With respect to this decomposition the derivation $A \in \Dera (\fg)$
is given by
$$ A= \left( \begin{array}{cccc}
0&0&0&1\\
-1 &0&0&0\\
0&0& U_1&0\\
0&0&0&0
             \end{array}
\right)  $$
where $\Span\{A_0,U_1\} \in \fso (\fg_0)$ is a 2-dimensional abelian subalgebra.\\
\end{enumerate}  
\end{pr}
We denote the corresponding simply connected Lie group with Lie algebra $\fd_A(\fg,\RR)$ by
$A(1,n-3)$ in case 1., by $\Osc(A_0,U_1)$ in case 2., and by $D(A_0,U_1)$ in case 3.

Proposition \ref{P7.1} is more a structure theorem than a full classification. We will give now an
exact classification of indecomposable metric Lie algebras of signature $(2,n-2)$ which have a
one-dimensional centre. These algebras are exactly those which appear in Proposition \ref{P7.1},
case 1. Therefore Proposition \ref{P6} implies

\begin{pr}\label{P7.2}
If $(\fd,\ip)$ is an indecomposable metric Lie algebra of signature $(2,n-2)$ with one-dimensional
centre. Then $\fd$ is isomorphic to
\begin{enumerate}
\item $\underline L_2(1,1)$ if $n=4$,
\item $\underline L_3(1,2)$ if $n=5$,
\item $\underline L_{2,\underline\lambda}(1,n-3)$ for $n>5$ even, or
\item $\underline L_{3,\underline\lambda}(1,n-3)$ for $n>5$ odd,
\end{enumerate}
where
\vspace{-2ex}
\begin{eqnarray*}
L_2:=\left( \begin{array}{cc}
0& 1 \\
1& 0
\end{array} \right)\,, 
\,L_3:= \left( \begin{array}{ccc}
0& 1 & 0\\
1&0& 1\\
0&-1&0
\end{array} \right) \,,\,
 L_{2,\underline \lambda}:= \left( \begin{array}{cc}
L_2& 0 \\
0& A_{\underline \lambda}
\end{array} \right)\,,\,
L_{3,\underline \lambda}:= \left( \begin{array}{cc}
L_3& 0 \\
0& A_{\underline \lambda}
\end{array} \right)\,,
\end{eqnarray*}

%
$\underline \lambda=(\lambda_1,\lambda_2,\dots,\lambda_l)$
with  
$0<\lambda_1\le\dots\le\lambda_l$ and $A_{\underline\lambda}$ is defined as in 
(\ref{EAlambda}).
\end{pr}
Using Theorem \ref{T5} we can compute the number of linearly independent parallel spinors on the
simply connected groups which correspond to the Lie algebras described in Proposition \ref{P7.2}.
The 4-dimensional Lie group $L_2(1,1)$ has 2 linearely independent parallel 
spinors. The 5-dimensional Lie group $L_3(1,2)$ of signature (2,3) 
is non-flat, Ricci-flat and has 3 linearely independent parallel spinors. The groups 
$L_{2,\underline\lambda}(1,n-3)$ and  $L_{3,\underline\lambda}(1,n-3)$, $n>5$, have 
$2^{\big[\frac n2\big]-1}$ linearly independent parallel spinors.\\[0.2cm]
Next we determine the space of parallel spinors on simply-connected Lie groups with Lie algebra
occuring in Proposition \ref{P7.1}, case 3. These Lie algebras are exactly those indecomposable
ones of signature $(2,n-2)$ which are odd-dimensional and have a 2-dimensional centre.

\begin{pr}\label{Ex5.5}
If $(\fd,\ip)$ is an odd-dimensional indecomposable metric Lie algebra of signature $(2,n-2)$ with
2-dimensional centre, then $n\ge 7$ and the simply connected Lie group with Lie algebra $\fd$ admits
$2^{\frac{n-5}2}$ linearly independent parallel spinors.
\end{pr}
\proof
Under the above assumptions $\fd$ is a double extension $\fd=\fd_A(\fg,\RR)$
where $(\fg,\ip_\fg)$ and $A\in \Dera(\fg,\ip)$ are as in Proposition \ref{P7.1}, 3. 
In particular,
\[\fg=\RR\oplus\fd_{A_0}(\fg_0,\RR)=\RR\oplus\RR\alpha_1\oplus\fg_0\oplus\RR H_1,\]
where $\fg_0$ is a $2m$-dimensional Euclidean and abelian Lie algebra and $A_0\in\fso(\fg_0)$ is 
bijective. $A$ is given by $U_1\in\fso(\fg_0)$. 
Obviously,
\[ 
[\fg,\fg]_\fg=[\fd_{A_0}(\fg_0,\RR),\fd_{A_0}(\fg_0,\RR)]_\fg=\RR\alpha_1\oplus{\rm 
Im} A_0=\RR\alpha_1\oplus\fg_0\]
and 
\[A\fg=\RR \oplus \RR\alpha_1 \oplus 
U_1(\fg_0)\,.\]
Therefore
\begin{equation}\label{Egg+Ag}
[\fg,\fg]_\fg+A\fg=\RR\oplus\RR\alpha_1\oplus\fg_0\,. 
\end{equation}
Now let $y$ be a unit vector in $\RR$ and $X_1,\dots,X_{2m}$ an 
orthonormal basis of $\fg_0$. By Proposition \ref{P5} and (\ref{Egg+Ag})
\begin{eqnarray*}
\widetilde{\hol}(D_A(\fg))&=&\Big\{\ \, 
2AZ\cdot\alpha+y\cdot[Z,y]_\fg+\sum_{j=1}^{2m}X_j\cdot[Z,X_j]_\fg\\
&&\quad-\frac 12(H_1-\alpha_1)[Z,H_1-\alpha_1]_\fg
+\frac 12(H_1+\alpha_1)[Z,H_1+\alpha_1]_\fg\quad \Big|\\
&& \hspace{8cm}Z\in [\fg,\fg]_\fg+A\fg\ \,\Big\}\\
&=&\Big\{\ \, 
2AZ\cdot\alpha+\sum_{j=1}^{2m}X_j\cdot[Z,X_j]_\fg+\alpha_1\cdot[Z,H_1]_\fg\ \Big| \ 
Z\in\RR y\oplus\RR\alpha_1\oplus\fg_0\ \Big\}\\
&=&\Span \,\{\ \alpha_1\cdot\alpha,\ 
U_1X\cdot\alpha+A_0X\cdot\alpha_1\mid X\in\fg_0 \}\,.
\end{eqnarray*}
Consider now $u=u_{++}\otimes u(1,1) + u_{+-}\otimes u(1,-1) + u_{-+}\otimes u(-1,1) 
+ u_{--}\otimes u(-1,-1) \in 
\Delta_{2m+1}\otimes\Delta_{2,2}=\Delta_{\RR y\oplus\fg_0}\otimes\Delta_{2,2}$. We have
\begin{eqnarray*}
\alpha_1\cdot\alpha\cdot u&=& -2u_{++}\otimes u(-1,-1)
\end{eqnarray*}
and
\begin{eqnarray*}
(U_1X\cdot\alpha+A_0X\cdot\alpha_1)\cdot u&=&\sqrt 2 \,\Big( -(U_1Xu_{++})\otimes 
u(1,-1)+(U_1Xu_{-+})\otimes u(-1,-1)  \\
&& \qquad+(A_0Xu_{++})\otimes u(-1,1)+ (A_0Xu_{+-})\otimes u(-1,-1)\,\Big).
\end{eqnarray*}
Therefore $\widetilde{\hol}(D_A(\fg))u=0\,$ is equivalent to
\begin{eqnarray}u_{++}=0\,,&
&U_1X\cdot u_{-+}+A_0X\cdot u_{+-}=0 \label{ea3}
\end{eqnarray}
for all $X\in\fg_0$.
Since $A_0$ and $U_1$ commute we may assume
\[\begin{array}{ll}
A_0X_{2i-1}=\lambda_iX_{2i},&A_0X_{2i}=-\lambda_iX_{2i-1}\\[1ex]
U_1X_{2i-1}=\mu_iX_{2i},&U_1X_{2i}=-\mu_iX_{2i-1}
\end{array}
\]
for $i=1,\dots,2m$. Therefore we obtain from (\ref{ea3})
\[\mu_iX_{2i}\cdot u_{-+}+\lambda_iX_{2i}u_{+-}=0\]
and, thus,
\[ \mu_iu_{-+}=-\lambda_iu_{+-}\,.\]
Since $U_1$ and $A_0$ 
are linearly independent this implies $u_{-+}=u_{+-}=0$.
Hence, 
\[u_{++}=u_{+-}=u_{-+}=0\,.\] 
We conclude
\[
\dim {\cal P}_{D_A}= \dim \Delta_{2m+1} =2^{m}  . 
\]
\qed

It remains to consider case 2.~of Proposition \ref{P7.1}. The Lie algebras occuring there are
exactly the indecomposable ones of signature $(2,n-2)$ which are even-dimensional and have a
2-dimensional centre. Since these spaces will occur in a more general context in Section
\ref{SSmult} we will postpone the calculation of the space of parallel spinors on $Osc (A_0, U_1)$.
However, we alredy mention the result in the following table.
\vspace{2ex}

{\bf Table 2: Indecomposable simply connected Lie groups with biinvariant metric of signature
$(2,n-2)$ and number of parallel spinors}
\begin{center}
\begin{tabular}{|l|c|c|}
\hline &&\\[-0.5ex]
$G$ & $\dim G =n$ & $\dim {\cal P}_G$  \\[2ex] \hline \hline &&\\[-0.5ex]
$(\widetilde{\SL} (2, \RR), -c B_{\fsl (2,\RR)}),\ c \in \RR^+$ & $3$ & $0$\\[2ex] 
\hline && \\[-0.5ex]
$L_2 (1,1)$ &$4$ & $2$ \\ [2ex]
\hline && \\[-0.5ex]
$L_3 (1,2) $ & $5$ & $3$\\[2ex] 
\hline && \\[-0.5ex]
$L_{2,\underline\lambda}(1,n-3) ,\ n>5,\ n \mbox{ even}$ & $n >5$ & 
$2^{\big[\frac{n}{2}\big]-1}$\\[2ex] 
\hline && \\[-0.5ex]
$L_{3,\underline\lambda}(1,n-3) ,\ n>5,\ n \mbox{ odd}$ & $n >5$ & 
$2^{\big[\frac{n}{2}\big]-1}$\\[2ex] 
\hline && \\[-0.5ex]
$Osc (A_0, U_1),\ A_0,U_1\in\fso(2m)$ & $2m+4 \ge 8$ & $2^m$ \\
$A_0,\,U_1$ as in Proposition \ref{P7.1}&&\\[2ex] 
\hline && \\[-0.5ex]
$D(A_0, U_1),\ A_0,U_1\in\fso(2m)$ & $2m+5 $ & $2^m\ $ \\
$A_0$ bijective; $A_0,\,U_1$ linearly independent&&\\[2ex] 
\hline
\end{tabular}
\end{center}

%
   \section{Indecomposable simply connected metric Lie groups 
            of dimension $n \le 6$}
%

Now we will classify the indecomposable simply connected metric Lie groups 
of dimension $n \le 6$ and determine the space of parallel spinors for each isomorphism class. We
will summarize the results in Table 3 - Table 6. Since we already have a classification for index 1
and 2 it remains to study the case of signature $(3,3)$. The Lie algebra $\fspin(1,3)$ is the only
simple Lie algebra, whose Killing form has this signature. Each indecomposable metric Lie algebra
$\fd$ which is not simple is a double extension by a simple or a one-dimensional Lie algebra $\fh$.
If $\fh$ is simple and $n\le 6$ , then $\fh$ must be 3-dimensional, hence, $\fh=\mathfrak{su}(2)$
or $\fh=\mathfrak{sl}(2,\RR)$. Consider now the case where $\fh$ is one-dimensional, i.e.
$\fd=\fd_A(\fg,\RR)$. Then $\fg$ is of
signature $(2,2)$. Since $\fd$ is indecomposable, $\fg$ cannot contain a simple ideal.
Consequently, $\fg=\underline L_2(1,1)$ or $\fg$ is abelian. It is not hard to prove that any
double extension of $\fg=\underline L_2(1,1)$ is decomposable (see also Remark \ref{R8}). This
implies $\fd=\fd_A(\fg,\RR)$, where $\fg$ is abelian and of signature $(2,2)$ and 
the kernel of $A \in \fso(2,2)$ is totally isotropic. Such a linear map $A$
is up to conjugation and up to non-zero multiples 
one of the following ones
\begin{enumerate}
\item $\quad A=N_1=\left({\small \begin{array}{cccc}
\textstyle 0&1&0&0\\0&0&0&0\\0&0&0&1\\0&0&0&0\end{array}} \right)\,,\quad 
\ip_\fg=\left({\small\begin{array}{cccc}
0&0&0&1\\0&0&-1&0\\0&-1&0&0\\1&0&0&0\end{array}} \right)$,
\item $\quad A=N_{2,t}=\left({\small\begin{array}{cccc}
1&-t&0&0\\t&1&0&0\\0&0&-1&-t\\0&0&t&-1\end{array} }\right)\,,\quad 
\ip_\fg=\left({\small\begin{array}{cccc}0&0&0&1\\0&0&-1&0\\0&-1&0&0\\1&0&0&0
\end{array}} \right)$,\ $\ t>0$,
\item $\quad A=N_{3,\pm}=\left({\small\begin{array}{cccc}
0&-1&\pm 1&0\\1&0&0&\pm 1\\0&0&0&-1\\0&0&1&0\end{array}} \right)\,,\quad 
\ip_\fg=\left({\small\begin{array}{cccc}
0&0&0&1\\0&0&-1&0\\0&-1&0&0\\1&0&0&0\end{array} }\right)$,
\item $\quad A=N_{4,t}=\left({\small\begin{array}{cccc}
0&-1&0&0\\1&0&0&0\\0&0&0&-t\\0&0&t&0\end{array} }\right)\,,\quad 
\ip_\fg=\left({\small\begin{array}{cccc}
-1&0&0&0\\0&-1&0&0\\0&0&1&0\\0&0&0&1\end{array}} \right)$, $\ t>0$,
\item $\quad A=N_5=\left({\small\begin{array}{cccc}
1&1&0&0\\0&1&0&0\\0&0&-1&1\\0&0&0&-1\end{array} }\right)\,,\quad 
\ip_\fg=\left({\small\begin{array}{cccc}
0&0&0&1\\0&0&-1&0\\0&-1&0&0\\1&0&0&0\end{array}} \right)$,
\item $\quad A=N_{6,t}=\left({\small\begin{array}{cccc}
1&0&0&0\\0&-1&0&0\\0&0&t&0\\0&0&0&-t\end{array}} \right)\,,\quad 
\ip_\fg=\left({\small\begin{array}{cccc}
0&1&0&0\\1&0&0&0\\0&0&0&1\\0&0&1&0\end{array}} \right)$,\ $\ t\ge 1$.
\end{enumerate}
Only $N_1$ has a non-trivial kernel.
Since $N_1^2=0$ the 6-dimensional indecomposable Lie group $N_1(2,2)$ of 
signature $(3,3)$ is flat and has 8 linearly independent parallel spinors. 
The 6-dimensional Lie groups $N_{k}(2,2),\ k=2,\dots,6$ have 4 linearly independent 
parallel spinors. $N_{k}(2,2)$ is Ricci-flat if and only if $k=2$ with the parameter $t=1$. 

\newpage

{\bf Table 3: $\dim G=3$} 
\begin{center}
\begin{tabular}{|l|c|c|l|}
\hline &&&\\[-0.5ex]
\parbox{5.7cm}{$G$} & signature  & $\dim {\cal P}_G$ 
&\parbox{5.3cm}{geometry}\\&$(p,q)$&& \\[2ex] \hline\hline &&&\\[-0.5ex]
 $(SU(2), -c B_{\fsu (2)}),\ c \in \RR^+$ & $(0,3)$ & 0& $G$ simple\\&&&simple 
holonomy \\[-1ex]\cline{1-3}&&&\\[-0.9ex]
 $(\widetilde{\SL} (2, \RR) , -c B_{\fsl (2, \RR)}),\ c \in \RR^+$ & $(1,2)$ & $0$ & 
 Einstein space, $R>0$\\[2ex] \hline
\end{tabular}
\end{center}
\vspace{4ex}
{\bf Table 4: $\dim G=4$}
\begin{center}
\begin{tabular}{|l|c|c|l|}
\hline &&&\\[-0.5ex]
\parbox{5.7cm}{$G$}  & signature  & $\dim {\cal P}_G$ 
&\parbox{5.3cm}{geometry}\\&$(p,q)$&& \\[2ex]\hline \hline &&&\\[-0.5ex]
 $Osc(1)$, i.e. the oscillator group  & $(1,3)$ & 2& $G$ solvable 
\\$Osc(\underline\lambda)$ with $\underline\lambda=\lambda_1=1$&&& abelian holonomy 
$\RR^2\subset\fso(1,3)$ \\&&&$R=0$,\ 
$\Ric(H,H)=-2$\\[2ex] \hline &&&\\[-0.5ex]
 $L_2(1,1)$ & $(2,2)$ & $2$ &  $G$ solvable \\ &&& 
abelian holonomy $\RR^2\subset\fso(2,2)$ \\ &&& $R=0$,\ $\Ric(H,H)=2$\\[2ex] 
\hline
\end{tabular}
\end{center}
\vspace{4ex}
   {\bf Table 5: $\dim G=5$}
\begin{center}
\begin{tabular}{|l|c|c|l|}
\hline &&&\\[-0.5ex]
\parbox{5.7cm}{$G$} & signature  & $\dim {\cal P}_G$ 
&\parbox{5.3cm}{geometry}\\&$(p,q)$&& \\[2ex] \hline\hline &&&\\[-0.5ex]
 $L_3(1,2)$ & $(2,3)$ & 3& $G$ 3-step nilpotent\\
&&&abelian holonomy $\RR\subset\fso(2,3)$ \\ &&&$\Ric=0$, non-flat\\[2ex] \hline\end{tabular}
\end{center}
\newpage
{\bf Table 6: $\dim G=6$}
\begin{center}
\begin{tabular}{|l|c|c|l|}
\hline &&&\\[-0.5ex]
$G$ & signature  & $\dim {\cal P}_G$ &geometry\\&$(p,q)$&& \\[2ex] \hline\hline &&&\\[-0.5ex]
$Osc(1,\lambda),\ \lambda\in\RR,\ \lambda\ge1$&$(1,5)$&4& $G$ solvable\\
&&&abelian holonomy $\RR^4$\\
&&&$\Ric(H,H)=-1-\lambda^2$\\[2ex] \hline &&&\\[-0.5ex]
$L_{2,\lambda}(1,3),\ \lambda\in\RR,\ \lambda>0$&$(2,4)$&4& $G$ solvable\\
&&&abelian holonomy $\RR^4$\\ 
&&&$\Ric(H,H)=2-2\lambda^2$\\[2ex] \hline &&&\\[-0.5ex]
$(Spin(1,3),-cB_{\fspin(1,3)})\,,\, c \in \RR^+$ & (3,3) & 0 & $G$ simple\\
&&& simple holonomy \\
&&& Einstein space, $R>0$ \\[2ex]\hline &&&\\[-0.5ex]
$T^*SU(2)_c,\ c\in\RR$ & $(3,3)$ & 1 & $G$ proper Levi\\
&&&holonomy equals $\fg$\\
[-1ex]\cline{1-3}&&& \\[-0.9ex]
$T^*\widetilde{SL}(2,\RR)_c,\ c\in\RR$ & $(3,3)$ & 1 &$\Ric^2=0,\ \Ric\not=0,\ R=0$\\[2ex]\hline
&&&\\[-0.5ex]
$N_1(2,2)$ &$(3,3)$ &8 & $G$ 2-step nilpotent\\
&&&holonomy $=0$, flat\\[2ex]\hline &&&\\[-0.5ex]
$N_k(2,2), k=2,\dots,6$ &$(3,3)$ &4 & $G$ solvable \\
 &&& abelian holonomy $\RR^4$ \\
&&& $\Ric^2=0, R=0$ \\
&&& $\Ric=0 \Leftrightarrow k=2,t=1$  \\[2ex] \hline 
\end{tabular}
\end{center}

\section
{Solvable Lie groups with maximal isotropic centre} 

\label{SSmult}
In this chapter we  study a special class of solvable metric Lie algebras with 
maximal isotropic centre. For that reason we consider  
multiply 1-dimensional extended Lie algebras.
Let $(\fg_0 , \ip_0)$ be an $n$-dimensional Euclidean abelian Lie algebra and $A_0 
=U_0 \in \fso (\fg_0)$ a bijective antisymmetric map.
Then the oscillator algebra
$$ \fg_1 := \fd_{A_0} (\fg_0 , \RR) $$
is an indecomposable Lorentzian Lie algebra with maximal isotropic centre 
$\fz(\fg_1)=\RR\alpha_1$. We now define a series of
Lie algebras $\fg_1 , \fg_2 , \fg_3 , \ldots$ by induction. Let $\fg_1 , 
\ldots , \fg_{m-1}$ be already defined. Choose $A_{m-1} \in \Dera (
\fg_{m-1} , \ip_{m-1})$ and define
\begin{eqnarray*}
\fg_m &:=& \fd_{A_{m-1}} (\fg_{m-1} ,\RR) \,=\,\RR \alpha_m\oplus\fg_{m-1}\oplus\RR 
H_m\\
&=& {\rm span} \{\alpha_1 , \ldots , \alpha_m , \fg_0 , H_1 , \ldots , H_m\}
\,. 
\end{eqnarray*}
We will denote the Lie bracket of $\fg_m$ by $[\,\cdot\,,\cdot\,]_m$.
If $\fg_m$ arises from $\fg_0$ by this procedure using $A_0,\dots, A_{m-1}$ we will 
write
\[ \fg_m=\fd(\fg_0,A_0,\dots,A_{m-1}).\]
The Lie algebra $\fg_m$ has a pseudo-Euclidean metric of signature $(m, m+n)$ given 
by
$$ \ip_m = \left( \begin{array}{ccc}
0&0 & E\\
0&\mbox{} \  \ip_0 \ \mbox{} &0\\
E&0&0 
\end{array}\right)$$
with respect to $\alpha_1 , \ldots , \alpha_m , \fg_0 , H_1 , \ldots , H_m
$. 

Consider now linear maps $A^0_k\in \fso(k,k+n),\ k=1,\dots,m-1$ and a bijective map 
$A_0^0\in\fso(n)$. Then $\{A_0^0,A_1^0,\dots,A_{m-1}^0\}$ is called a {\em normal 
set of derivations} if
\begin{itemize}
\item[(i)]
$A_k^0 = \left( \begin{array}{ccc}
0&0&Z_k\\
0&U_k&0\\
0&0&0 \end{array} \right) 
\,$, 
where $U_k \in \fso (n), Z_k \in \RR(k)$ with $Z^{\top}_k = - Z_k$ and
\item[(ii)]
${\rm span} \{U_0, \ldots , U_{m-1}\} \subset \fso (n)$ is abelian.
\end{itemize} 
If as above $\fg_0$ is a Euclidean abelian Lie algebra, then one proves by induction 
that $A^0_k$ is indeed a derivation on $\fg_k^0=\fd(\fg_0, A_0^0,\dots,A_{k-1}^0)$. 


\begin{pr}\label{P6.1}
Let $\{A_0^0, \ldots ,A^0_{m-1}\}$ be a set of normal derivations. Then the Lie 
bracket of 
\begin{eqnarray*}
\fg_m^0 &=& \fd(\fg_0,A_0^0,\dots,A_{m-1}^0)\\
&=&\RR \alpha_m \oplus \ldots \oplus \RR \alpha_1
\oplus \fg_0 \oplus \RR H_1 \oplus \ldots \oplus \RR H_m
\end{eqnarray*}
is given by
\begin{eqnarray}\label{E12}
[\alpha_i , \cdot\, ]_m &=& 0\nonumber \\[1ex]
\ [X,Y]_m &=& \sum\limits^m_{j=1} \langle U_{j-1} X , Y \rangle_0 \alpha_j
\quad \quad \quad X,Y \in \fg_0 \nonumber \\[1ex]
[H_i , X]_m &=& U_{i-1} X \quad \quad \quad i=1, \ldots , m\\[1ex]
[H_i , H_j]_m &=& (Z_{m-1})^j_i \alpha_{m} + (Z_{m-2})^j_i
\alpha_{m-1} + \ldots \nonumber \\
&& + (Z_i)^j_i \alpha_{i+1} + \sum\limits^{i-1}_{k=1} (Z_{i-1})^k_j
\alpha_k \quad \quad \mbox{{\rm for} $m \ge i >j \ge 1$} \ . \nonumber 
\end{eqnarray}
In particular, $\fg_m^0$ is solvable. Let $G^0_m$ be the simply connected metric Lie 
group associated with $(\fg_m^0,\ip_m)$. Then
\begin{equation}\label{Eholg0} \hol (G^0_m)=  \left\{ \ \left( \begin{array}{ccc}
0& \begin{array}{c}(U_{m-1}X)^{\top}\\\vdots\\(U_0X)^{\top} \end{array}
&0\\[2ex]
0&0& -U_0X \ldots -U_{m-1}X\\[1em]
0&0&0 \end{array} \right) \quad \quad \Bigg| \quad X \in \fg_0\  \right\}
\end{equation}
is abelian and for the spin holonomy
\begin{equation}\label{spih}
\widetilde{\hol} (G^0_m) = \Big\{\ \sum\limits^m_{k=1} (U_{k-1} X) \cdot
\alpha_k\mid X\in\fg_0\ \Big\} \subset \fspin (\fg_m^0 )
\end{equation}
holds. The Ricci tensor of $G_m^0$  is given by
$$\Ric (H_i , H_j)= {\tr} (U_{i-1} \cdot U_{j-1}) \ . $$
All other components are zero. In particular, $G_m^0$ is a scalar flat,
non-Einstein space with 2-step nilpotent Ricci tensor. 
\end{pr}
\proof 
Relation (\ref{E12}) follows from the definition of $\fg_m^0$; the
formula for the Ricci tensor from Proposition \ref{P2}. Obviously, the third 
derivative
of $\fg_m^0$ vanishes, hence $\fg_m^0$ is solvable.
The assertion on $\hol (G^0_m)$ follows from
\[ \hol (G^0_m)= 
\ad_m([\fg^0_m,\fg^0_m])=\ad_m(\Span\{\alpha_1,\dots,\alpha_m\}\oplus 
\fg_0)=\ad_m(\fg_0)\,.\]
and (\ref{E12}). Equation (\ref{spih}) is a consequence of (\ref{E8}) and (\ref{Eholg0}).
\qed
\ \\
If $\{A_0^0,\dots, A_{m-1}^0\}$ is a normal set of derivations such that \\[0.1cm]
$$ A^0_k = \left(\begin{array}{ccc}0&0&0\\0&U_{k}&0\\0&0&0\end{array}\right)\;, 
\;k=0,\dots, m-1 $$\\[0.1cm]
we call the Lie algebra $\fg_m^0=\fd(\fg_0,A^0_0,\dots,A_{m-1}^0)$ a {\em 
generalized oscillator algebra of index m} and denote it by 
\[ \osc(U_0,\dots,U_{m-1}) .\]
The corresponding simply-connected Lie group is denoted by 
$\Osc(U_0,\dots,U_{m-1})$.
\begin{co}
For a set $\{A_0^0, \ldots ,A^0_{m-1}\}$ of normal derivations we may assume that 
all $U_k$ with $U_k \not= 0$ are linearly independent.
\end{co}
\proof
Fix $k<m$ and assume that all $U_l$, $l<k$, with $U_l\not=0$ are linearly 
independent. If $U_k$ is a linear combination of these $U_l$, i.e. if
\[ U_{k}=\lambda_0 U_0+\dots+\lambda_{k-1}U_{k-1},\qquad 
\lambda_0,\dots,\lambda_{k-1}\in\RR\,,\]
then we consider the inner derivation
\[\ad_{k}\Big(\sum_{j=1}^{k}\lambda_{j-1}H_j\Big)=\sum_{j=1}^{k}\lambda_{j-1}
\left(\begin{array}{ccc}0&0&Z_{j-1}\\0&U_{j-1}&0\\0&0&0\end{array}\right)
=\left(\begin{array}{ccc}0&0&\tilde Z_{k}\\0&U_{k}&0\\0&0&0\end{array}\right)\]
and therefore
\[\left(\begin{array}{ccc}0&0&Z_{k}\\0&U_{k}&0\\0&0&0\end{array}\right)\equiv
\left(\begin{array}{ccc}0&0&Z_k-\tilde 
Z_{k}\\0&0&0\\0&0&0\end{array}\right)\qquad{\rm mod}\ \ad(\fg_{k})\,.\]
The assertion now follows from Proposition \ref{P6}.
\qed

We know from Proposition \ref{P6.1} that the Lie algebra $\fg_m^0$ constructed above 
is solvable and $\Span\{\alpha_1,\dots,\alpha_m\}\subset \fz(\fg_m^0)$ is a maximal 
isotropic subspace of $\fg_m^0$. Conversely, the following is true 


\begin{pr} \label{P6.2}
Let $(\fg_m,\ip)$ be a solvable metric Lie algebra of signature $(m,m+n)$, 
$m,n\in\NN$, such that the centre $\fz(\fg_m)$ of $\fg_m$ contains a maximal 
isotropic subspace $\fa_m$. Then
\begin{enumerate}
\item
$(\fg_m,\ip)$ is isomorphic to a Lie algebra $\fd(\fg_0,A_0,\dots,A_{m-1})$, where 
$\fg_0$ is a Euclidean abelian Lie algebra, $\fa_m=\Span\{\alpha_1,\dots,\alpha_m\}$ 
and the derivations $A_k$ satisfy
$\, A_k(\Span\{\alpha_1,\dots,\alpha_k\})=0 \;$ for $\; k=1,\ldots, m-1$.
\item
If $A_0$ in {\it 1.} is bijective, then 
$\,\fg_m = \fd(\fg_0,A_0^0,\dots,A_{m-1}^0)\,$
for a normal set of derivations $\{A_0^0,\dots,A_{m-1}^0\}$.
\end{enumerate}
\end{pr}

\proof
Let $\fa_m \subset \fz(\fg_m)$ be maximal isotropic. Choose an element $\alpha_m \in 
\fa_m$ with $\alpha_m \not= 0$.
Then there exists an element $H_m\in \fg_m$ such that $\langle\alpha_m,H_m\rangle=1$ 
and $\langle H_m,H_m\rangle=0$. Define $\fg_{m-1}:=\{\alpha_m,H_m\}^\perp$. Then 
there is a Lie algebra structure on $\fg_{m-1}$ such that
$$ \fg_m = \fd_{A_{m-1}} (\fg_{m-1} , \RR)= \RR \alpha_m \oplus
\fg_{m-1} \oplus \RR H_{m}\,, $$ 
where $A_{m-1}=\ad_m(H_m)|_{\fg_{m-1}}$.
Since $\fg_m$ is solvable, $\fg_{m-1}$ is solvable as well.

By (\ref{Ec}) the center of $\fg_{m}$ equals
\begin{equation}\label{Ezg} \fz (\fg)= \left\{ \begin{array}{l}
\RR \alpha_m \oplus (\ker A_{m-1} \cap \fz (\fg_{m-1})) \ \mbox{ if }  \ 
A_{m-1}\not\equiv 0 \mbox{ mod } \ad(\fg_{m-1})\\ 
\RR \alpha_m \oplus  \fz (\fg_{m-1})\oplus \RR(H_m-X_0) \ \mbox{ if } \ 
A_{m-1}=\ad_{m-1}(X_0)\ .
\end{array}\right. \end{equation}
Since $\fa_m \in \fz(\fg_m)$ is totally isotropic and $\alpha_m\in\fa_m$ we have 
$\fa_m\subset\RR\alpha_m\oplus\fg_{m-1}$. Hence, we obtain 
$\fa_m=\RR\alpha_m\oplus\fa_{m-1}$, where $\fa_{m-1}\subset\fg_{m-1}$ is maximal 
isotropic.
On the other hand, it follows from $\fa_m\subset\fz(\fg)$ and (\ref{Ezg}) that
\[\fa_{m-1}\subset\fz(\fg_{m-1})\cap\Ker A_{m-1} \subset \fz(\fg)\,.\]
Hence the center $\fz(\fg_{m-1})$ admits a maximal isotropic subspace $\fa_{m-1}$. 
If we continue with $\fa_{m-1}$ in the same way as with $\fa_m$ we obtain a sequence 
of isotropic elements $\alpha_m,\alpha_{m-1},\dots,\alpha_1 \in \fa_m$ and 
derivations $A_{m-1},A_{m-2},\dots,A_1$ such that 
\[
\fg_m=\fd(\fg_0,A_0,\dots,A_{m-1}) \hspace{0.5cm}\mbox{and}\hspace{0.5cm} 
A_k(\Span\{\alpha_1,\dots,\alpha_k\})=0 \;\;\; k=1,\dots,m-1.
\]
Now, let us assume in addition that the derivation $A_0$ in this representation is 
bijective. Then, if $m=1$, the second statement of the Proposition is obviously 
true. Let us suppose, that statement 2. is true for $m-1\geq 1$. Then we may assume 
that $\fg_m=\fd_{A_{m-1}}(\fg_{m-1},\RR)$ and
\begin{equation}
\label{Egm-1} 
\fg_{m-1} = \fd(\fg_0,A_0^0,\dots,A_{m-2}^0) 
\end{equation} 
where $\fg_0$ is a Euclidean abelian Lie algebra, $\{A_0^0,\dots,A_{m-2}^0\}$ is a 
normal set of derivations and $\fa_{m-1}=\Span\{\alpha_1,\dots,\alpha_{m-1}\}$. From 
$\fa_{m-1}\subset\Ker A_{m-1}$ we get
\begin{equation}\label{EA}
A_{m-1} = \left( \begin{array}{ccc}
0& -V^{\top}_{m-1} & Z_{m-1}\\
0&U_{m-1} & V_{m-1} \\
0&0&0 \end{array} \right) \ , \quad \begin{array}{l}
U_{m-1} \in \fso (n)\\
V_{m-1} \in \RR (n,m-1)\\
Z_{m-1} \in \RR (m-1), Z_{m-1}^{\top} =- Z_{m-1} \end{array} 
\end{equation}
with respect to $\alpha_1 , \ldots , \alpha_{m-1} , \fg_0 , H_1 , \ldots , H_{m-1}
$.
Let us prove, that there exists an antisymmetric derivation $A_{m-1}^0$ on 
$\fg_{m-1}$ such that $A_{m-1}$ is equal to $A_{m-1}^0$ modulo inner derivations and 
$\{A_0^0,\dots,A_{m-1}^0\}$ is a normal set of derivations. From (\ref{E12}) we know 
that
\begin{eqnarray*}
A_{m-1} [H_i , X]_{m-1}\, =\, A_{m-1} U_{i-1} X
\,= \,- \sum\limits^{m-1}_{j=1} \langle (V_{m-1})_j , U_{i-1} X \rangle_0 \alpha_j+
U_{m-1} U_{i-1} X
\end{eqnarray*}
and
\begin{eqnarray*}
\lefteqn{[A_{m-1} H_i, X]_{m-1} + [H_i , A_{m-1} X]_{m-1}=} \\
&=& 
 [\sum\limits^{m-1}_{l=1}(Z_{m-1})^l_{i} \alpha_l + (V_{m-1})_i , X]_{m-1} - 
\sum\limits^{m-1}_{l=1}\langle (V_{m-1})_l , X \rangle_0
[ H_i, \alpha_l]_{m-1} 
\\&&
+ [H_i , U_{m-1} X]_{m-1} \\
&=& \sum\limits^{m-1}_{j=1} \langle U_{j-1} (V_{m-1})_i, X \rangle_0 \alpha_j
+ U_{i-1} U_{m-1} X  
\end{eqnarray*}
for all $X\in \fg_0$.
It follows, that
$$ [U_j , U_{m-1}]=0 \quad \quad \mbox{for all $j=0, \ldots , m-2$} $$
and that
$$ U_{i-1} (V_{m-1})_j =U_{j-1} (V_{m-1})_i \quad \quad i,j =1 , \ldots, m-1 \ . $$
Hence, ${\Span} \{U_0 , \ldots , U_{m-1}\}$ is abelian and the matrix
$V_{m-1}$ is given by the columns
$$V_{m-1} = (Y, \ U^{-1}_0 U_1 Y, \ U^{-1}_0 U_2 Y, \ \ldots, \ U^{-1}_0 U_{m-2} Y)
\ , $$
where $Y=(V_{m-1})_1$.
By (\ref{E12}) we have
\begin{equation}\label{**} \ad_{m-1} (Y)= \left( \begin{array}{ccc}
0& \begin{array}{c} (U_{0}Y)^{\top}\\\vdots\\(U_{m-2}Y)^{\top} \end{array}
&0\\[2em]
0&0& -U_0Y \ldots - U_{m-2}Y\\[1em]
0&0&0 \end{array} \right) \,. 
\end{equation}
Therefore,
\begin{eqnarray*}
A_{m-1} &=& \left( \begin{array}{ccc}
0&0& Z_{m-1}\\
0& U_{m-1} &0\\
0&0&0 \end{array} \right) \,-\,\ad_{m-1} (U_0^{-1} (V_{m-1})_1)\\[1ex]
&\equiv& \left( \begin{array}{ccc} 
0&0& Z_{m-1}\\
0& U_{m-1} &0\\
0&0&0 \end{array} \right) \ {\rm mod} \ {\rm ad} (\fg_{m-1}) \ . 
\end{eqnarray*}
\qed

For the following Proposition we need a simple fact on the Lie algebra $\fso(n)$.

\begin{lm} \label{L2}
If $\fb \subset \fso (n)$ is an abelian subalgebra and $[X, \fb ] \subset
\fb$ for $X\in\fso(n)$, then $[X, \fb]=0$. 
\end{lm}
\proof
Since $\fso (n)$ admits a definite inner product also the subalgebra  
$\mathfrak{c}=\fb\oplus\RR X$ has a definite inner product. By our assumption $\fb$ 
is an ideal in $\mathfrak{c}$. Hence, also $\fb^\perp$ is an  ideal and 
$\mathfrak{c}=\fb\oplus\fb^\perp$. Thus $\fb^\perp=\RR(X+B_0)$ for some $B_0\in\fb$ 
and therefore $0=[\fb,X+B_0]=[\fb,X]$.\qed
\begin{pr}\label{P6.3}
Let $(\fg_0,\ip_0)$ be an $n$-dimensional Euclidean abelian Lie algebra and 
$\{A_0^0,\dots,$ $A_{m-1}^0\}$ a normal set of derivations such that $U_0 , \ldots , 
U_{m-1} \in \fso (n)$ are linearly independent. If $A_m$ is an antisymmetric 
derivation on $\fg^0_m=\fd(\fg_0,A^0_0,\dots,A^0_{m-1})$ then the centre of the 
solvable Lie algebra $\fg_{m+1}=\fd_{A_m}(\fg_m^0,\RR)$ has a maximal isotropic 
subspace, i.e. an isotropic subspace of dimension $m+1$.
In particular, $\fg_{m+1}$ satisfies the condition of Proposition \ref{P6.2}, 2.
\end{pr}
\proof
If $U_0 , \ldots , U_{m-1}$ are linearly independent, then the centre of $\fg^0_m$ 
is equal to $\Span\{\alpha_1,\dots,\alpha_m\}$. Indeed, if 
$X+\sum_{i=1}^m\lambda_iH_i\in\fz(\fg^0_m)$
for $X\in\fg_0$, $\lambda_i\in\RR$, then we obtain from (\ref{E12}) 
\[0=[X+\sum_{i=1}^m\lambda_iH_i,Y]=\sum_{j=1}^m\langle 
U_{j-1}X,Y\rangle_0\alpha_j+\sum_{i=1}^m\lambda_i U_{i-1}Y\]
for all $Y\in\fg_0$. Hence, $X\in \Ker U_{j-1}$ for $j=1,\dots,m$. In particular, 
$U_0X=0$ which implies $X=0$. Furthermore, 
\[\sum_{i=1}^m\lambda_i U_{i-1}Y=0\]
for all $Y\in\fg_0$. Since $U_0 , \ldots , U_{m-1}$ are linearly independent we 
obtain $\lambda_i=0$ for $i=1,\dots,m$.

Now let $A_m$ be an antisymmetric derivation on $\fg_m^0$. Then
$$ A_m = \left( \begin{array}{ccc}
B & -V^{\top} &Z_m\\
0& U_m & V\\
0&0&-B^{\top} \end{array} \right) \quad \begin{array}{l}
U_m \in \fso (n)\\ (Z_m)^{\top} = - Z_m \end{array} $$
since $A_m \fz (\fg_m^0) \subset \fz (\fg_m^0)$. Using $\fz (\fg_m^0) =
{\rm span} \{\alpha_1 , \ldots, \alpha_m\}$ we obtain for $X,Y \in \fg_0$
\begin{eqnarray*}
A_m [X,Y]_m &=& \sum\limits^m_{j=1} \langle U_{j-1} X,Y \rangle_0 A_m \alpha_j
\,=\, \sum\limits^m_{k,j=1} \langle U_{j-1} X,Y \rangle_0 B_j^k \cdot
\alpha_k
\\[2ex]
[A_m X,Y]_m + [X, A_m Y]_m &=& [U_m X,Y]_m + [X,U_m Y]_m\\
&=& \sum\limits^m_{k=1} (\langle U_{k-1} U_m X,Y\rangle_0+
\langle U_{k-1} X, U_m Y \rangle_0 ) \alpha_k\\
&=& \sum\limits^m_{k=1} \langle [U_{k-1} , U_m] X, Y \rangle_0
\alpha_k
\end{eqnarray*}
It follows, that for all $k=1 , \ldots , m$
\begin{equation}\label{E15}
\sum\limits^m_{j=1} B^k_j \cdot U_{j-1} = [U_{k-1} , U_m ]_{{\fso}(n)}\,.
\end{equation}
We know, that $\fb:={\rm span} \{U_0 , \ldots , U_{m-1}\} \subset \fso (n)$ is
abelian. By (\ref{E15}) $[U_m, \fb] \subset \fb$. Now Lemma \ref{L2} implies $[U_m , 
\fb]=0$.
Since we supposed, that $U_0, \ldots , U_{m-1}$ are
linearly independent in $\fso (n)$, we obtain $B  \equiv 0$ by (\ref{E15}).
Therefore, $A_m \fz (\fg_m^0)=0$ and, consequently, ${\rm span}\{\alpha_1, 
\ldots , \alpha_{m+1}\}\subset\fz (\fg_{m+1})$ is a maximal isotropic subspace. \qed
\begin{re} If $\fg_1 =\fd_{A_0} (\fg_0, \RR)$,
$\fg_2= \fd_{A_1} (\fg_1 , \RR)$, $\fg_3 =\fd_{A_2} (\fg_2 , \RR)$ is a
sequence of indecomposable Lie algebras, then
$\fg_1$, $\fg_2$ and $\fg_3$ have a maximal isotropic centre.  \end{re}

For $\fg_1$ this follows from the definition. For $\fg_2$ and $\fg_3$ this
follows from the indecomposability of $\fg_2$. Indeed, if we take $U_1 =
\lambda U_0$, then 
$A_1= \ad_1
(\lambda H_1)$. Therefore $A_1$ would be an inner derivation and $\fg_2$
decomposable (Proposition \ref{P1}). Hence $U_0$ and $U_1$ have to be
linearly independent.  Consequently, by  Proposition \ref{P6.3} the
centres $\fz(\fg_2)$ and $\fz(\fg_3)$ of $\fg_2$ and $\fg_3$, respectively, contain 
a maximal isotropic subspace. On the other hand, $\fz(\fg_2)$ and $\fz(\fg_3)$ are 
totally isotropic since $\fg_2$ and $\fg_3$ are indecomposable.
\qed

\begin{re}\label{R8}
The previous remark shows that a double extension of a 4-dimensional oscillator 
algebra by an 1-dimensional factor is decomposable. With the same argument one shows 
that the 6-dimensional group of split signature arising by double extension of 
$L_2(1,1)$ is decomposable as well. 
\end{re}
Now, we describe the space of parallel spinors for the solvable metric Lie algebras 
with maximal isotropic center of the form $\fg_m^0$.

\begin{theo}\label{T6} Let $(\fg_0,\ip_0)$ be an abelian Euclidean Lie algebra of 
dimension $n$ and consider
 $\fg=\fd(\fg_0, A_0^0,\dots,A_{m-1}^0)$ for normal set of derivations 
$\{A_0^0,\dots,A_{m-1}^0\}$. Then the simply connected Lie group $G$ with Lie 
algebra $\fg$ has
exactly $2^{K+\frac{n}{2}}$ linearly independent parallel spinor fields, where 
\[K=m-\dim(\Span(\{U_0,\dots,U_{m-1}\})\,.
\] 
\end{theo}
\proof We may assume that all $U_k$ with $U_k\not=0$ are linearly independent. 
Hence, $K=\#\{k\mid U_k=0\}$. The parallel spinor fields correspond to the spinors
$$ v= \sum\limits_{\gep =(\gep_m, \ldots , \gep_1)} v_{\gep} \otimes u(
\gep_m , \ldots , \gep_1) \in \Delta_n \otimes \Delta_{m,m}$$
satisfying
$$ \sum\limits^m_{k=1} (U_{k-1} X) \cdot \alpha_k \cdot v =0$$
for all $X\in\fg_0$.
The dimension $n$ of $\fg_0$ is even since $U_0$ bijective. Let $n=2N$. Since $\Span 
\{U_0,\dots,U_{m-1}\}$ is abelian there exists a basis $e_1,\dots,e_n$ such that
\[U_k=\left( \begin{array}{ccc}
\Lambda_1^{(k)} && 0\\
& \ddots &\\
0&& \Lambda_m^{(k)}
             \end{array}
\right) \ ,\ \mbox{ where } 
\Lambda_r^{(k)}=\left(\begin{array}{cc} 0&-\lambda_r^{(k)}\\ \lambda_r^{(k)} &0
\end{array}\right)\,\]
for each $k\in \{0,\dots,m-1\}$.
Hence,
\[\sum^m_{k=1}(U_{k-1}X)\cdot\alpha_k\Big( \sum_{\eps}v_\eps\otimes 
u(\eps_m,\dots,\eps_1)\Big)=0\]
is satisfied for all $X\in\fg_0$ if and only if
\[\sum_{k=1}^m\sum_{\eps}(-1)^{k-1}(\eps_k+1)\eps_k\cdot\dots\cdot\eps_m 
\lambda^{(k-1)}_j X_{2j}v_\eps\otimes u(\eps_m,\dots,-\eps_k,\dots,\eps_1)=0\]
for $j=1,\dots,N$ and an orthonormal basis $X_1,\dots,X_{2N}$ of $\fg_0$. This is 
equivalent to 
\[\sum_{k=1}^m(-1)^k(+\eps_k+1)\eps_k\cdot\dots\cdot\eps_m \lambda^{(k-1)}_j 
v_{(\eps_m,\dots,\eps_k,\dots,\eps_1)}=0\]
for $j=1,\dots,N$ and for all $\eps$.  Since $\{U_k\mid U_k\not=0\}$ are linearly 
independent this is equivalent to
\[U_{k-1}\not=0,\ \eps_k=1 \ \Rightarrow\  v_{(\eps_m,\dots,\eps_1)}=0\,.\]
Consequently, 
\[\dim V_{\mathfrak{hol} (G_m)}=(\dim \Delta_n)\cdot 2^K =2^{K + \frac{n}{2}}\,.\]
\qed

\begin{re} 
In particular, each generalized oscillator group $\Osc(U_0,\dots,U_{m-1})$ of index 
$m$, where $\Span\{U_0,\dots,U_{m-1}\} \subset \fso(2l)$ is an $m$-dimensional 
abelian subalgebra, has $2^l$ linearly independent parallel spinors.
\end{re}

{\bf Acknowledgement} We wish to thank D.~Alekseevski for his helpful comments on the first version
 of this paper (SFB 288 preprint Nr. 543, 2002). Furthermore, we thank M.~Olbrich for valuable discussions and for
pointing out an error in a previous version of the manuscript.

\bibliographystyle{alpha}

\begin{thebibliography}{BFOHP01}

\bibitem[Ast73]{Astrakhantsev:73}
V.V. Astrakhantsev.
\newblock Pseudo-riemannian symmetric spaces with commutative holonomy group.
\newblock {\em Mat. Sbornik (russ)}, 90(2):288--305, 1973.

\bibitem[B{\"a}r93]{Baer:93}
Ch. B{\"a}r.
\newblock Real {K}illing spinors and holonomy.
\newblock {\em Comm. Math. Phys.}, 154:509--521, 1993.

\bibitem[Bau89]{Baum2:89}
H.~Baum.
\newblock Complete {R}iemannian manifolds with imaginary {K}illing spinors.
\newblock {\em Ann. Glob. Anal. Geom.}, 7:205--226, 1989.

\bibitem[Bau00]{Baum1:00}
H.~Baum.
\newblock Twistor spinors on {L}orentzian symmetric spaces.
\newblock {\em Journ. Geom. Phys.}, 34:270--286, 2000.

\bibitem[BBI93]{Berard-Bergery/Ikemakhen:93}
L.~Berard-Bergery and A.~Ikemakhen.
\newblock On the holonomy of {L}orentzian manifolds.
\newblock {\em Proc. Symp. Pure Math.}, pages 27--40, 1993.

\bibitem[BFGK91]{Baum/Friedrich/ua:91}
H.~Baum, T.~Friedrich, R.~Grunewald, and I.~Kath.
\newblock {\em Twistors and {K}illing Spinors on {R}iemannian Manifolds},
  volume 124 of {\em Teubner-Texte zur Mathematik}.
\newblock Teubner-Verlag, Stuttgart/Leipzig, 1991.

\bibitem[BFOHP01]{Figueroa2:01}
M.~Blau, J.M. Figueroa-O'Farrill, C.~Hull, and G.~Papadopoulos.
\newblock New maximally supersymmetric background of {IIB} superstring theory.
\newblock hep-th/0110242, 2001.

\bibitem[Bie92]{Bieliavsky:92}
P.~Bieliavsky.
\newblock Extensions de dimension un de l'algebre de {H}eisenberg et espaces
  pseudoriemanniens symetrique.
\newblock {\em Acad. R. Belg.}, 3(12):299--315, 1992.

\bibitem[BK99]{Baum/Kath:99}
H.~Baum and I.~Kath.
\newblock Parallel spinors and holonomy groups on pseudo-{R}iemannian spin
  manifolds.
\newblock {\em Ann. Glob. Anal. Geom.}, 17:1--17, 1999.

\bibitem[Boh00]{Bohle:00}
Ch. Bohle.
\newblock Killing spinors on {L}orentzian manifolds.
\newblock {\em Journ. Geom. Phys.}, 2000.

\bibitem[Bor97]{Bordemann:97}
M.~Bordemann.
\newblock Nondegenerate invariant bilinear forms on nonassociative algebras.
\newblock {\em Acta Math.Univ. Comenianae}, LXVI(2):151--201, 1997.

\bibitem[Bou00]{Boubel:00}
Ch. Boubel.
\newblock Sur l'holonomie des varietes pseudo-riemanniennes.
\newblock {\em These, Univ. Nancy,}, 2000.

\bibitem[Bry00]{Bryant:00}
R.L. Bryant.
\newblock Pseudo--{R}iemannian metrics with parallel spinor fields and
  vanishing {R}icci tensor, in: {G}lobal {A}nalysis and {H}armonic {A}nalysis,
  eds. {J.P.B}ourguignon, {T.B}ranson, {O.H}ijazi.
\newblock {\em Seminaires et Congres, French Math. Soc.}, 4:53--94, 2000.

\bibitem[Buc00a]{Buchholz2:00}
V.~Buchholz.
\newblock A note on real {K}illing spinors in {W}eyl geometry.
\newblock {\em Journ. Geom. Phys.}, 35:93--98, 2000.

\bibitem[Buc00b]{Buchholz1:00}
V.~Buchholz.
\newblock Spinor equations in {W}eyl geometry.
\newblock {\em Suppl. di Rend. Circ. Mat. Palermo}, Ser.II, Nr.63:63--73, 2000.

\bibitem[CP80]{Cahen/Parker:80}
M.~Cahen and M.~Parker.
\newblock Pseudo-{R}iemannian symmetric spaces.
\newblock {\em Mem. AMS}, 24(229):1--108, 1980.

\bibitem[CW70]{Cahen/Wallach:70}
M.~Cahen and N.~Wallach.
\newblock Lorentzian symmetric spaces.
\newblock {\em Bull. AMS}, 76(3):585--591, 1970.

\bibitem[FI01]{Friedrich/Ivanov:01}
T.~Friedrich and S.~Ivanov.
\newblock Parallel spinors and connections with skew-symmetric torsion in
  string theory.
\newblock math.DG/0102142, 2001.

\bibitem[FO99a]{Figueroa:99}
J.M. Figueroa-O'Farrill.
\newblock Breaking the {M}-waves.
\newblock hep-th/9904124, 1999.

\bibitem[FO99b]{Figueroa2:99}
J.M. Figueroa-O'Farrill.
\newblock More {R}icci-flat branes.
\newblock ESI Preprint 769, 1999.

\bibitem[FOS95]{Figueroa:95}
J.M. Figueroa-O'Farrill and S.~Stanciu.
\newblock On the structure of symmetric self-dual {L}ie algebras.
\newblock hep-th/9506152, 1995.

\bibitem[FS87]{Favre:87}
G.~Favre and L.J. Santharoubane.
\newblock Symmetric, invariant, non-degenerate bilinear form on a {L}ie
  algebra.
\newblock {\em Journ. of Algebra}, 105:451--464, 1987.

\bibitem[Ike96]{Ikemakhen:96}
A.~Ikemakhen.
\newblock Examples of indecomposable non-irreducible {L}orentzian manifolds.
\newblock {\em Ann. Sci. Math. Quebec}, 20(1):53--66, 1996.

\bibitem[Ike99]{Ikemakhen1:99}
A.~Ikemakhen.
\newblock Sur l'holonomie des varietes pseudo-{R}iemannian de signature
  (2,2+n).
\newblock {\em Preprint, to appear in Publications Mathematiques}, 1999.

\bibitem[Kat99]{Kath:99}
I.~Kath.
\newblock Killing {S}pinors on {P}seudo-{R}iemannian {M}anifolds.
\newblock Habilitationsschrift {H}umboldt-{U}niversit{\"a}t {B}erlin, 1999.

\bibitem[Kat00]{Kath3:98}
I.~Kath.
\newblock Parallel {P}ure {S}pinors on {P}seudo-{R}iemannian {M}anifolds.
\newblock Sfb 288 preprint Nr. 356 (1998) and in: Topology and Geometry of
  Submanifolds X, 2000.

\bibitem[KO]{Kath/Olbrich:02}
I.~Kath and M.~Obrich.
\newblock Metric {L}ie algebras with maximal isotropic center.
\newblock In preparation.

\bibitem[Lei01]{Leistner:01}
Th. Leistner.
\newblock Lorentzian manifolds with special holonomy and parallel spinors.
\newblock to appear in Proceedings of the 21st Winter School on "Geometry and
  Physics" (Srni 2001), Rend. Circ. Mat. Palermo, 2001.

\bibitem[Med85]{Medina:85}
A.~Medina.
\newblock Groupes de {L}ie munis de metriques bi-invariantes.
\newblock {\em Tohoku Math. Journ.}, 37:405--421, 1985.

\bibitem[Mor96]{Moroianu:96}
A.~Moroianu.
\newblock Structures de {W}eyl admettant des spineurs parallelles.
\newblock {\em Bull. Soc. Math. France}, 1996.

\bibitem[MR85]{Medina/Revoy2:85}
A.~Medina and Ph. Revoy.
\newblock Algebres de {L}ie et produit scalaire invariant.
\newblock {\em Ann. scient.Ecole Norm. Sup., 4. serie}, 18:553--561, 1985.

\bibitem[MS00]{Moroianu/Semmelmann:00}
A.~Moroianu and U.~Semmelmann.
\newblock Parallel {S}pinors and holonomy groups.
\newblock {\em J. Math. Phys.}, 41:2395--2402, 2000.

\bibitem[Neu02]{Neukirchner:02}
Th. Neukirchner.
\newblock Pseudo-{R}iemannian symmetric spaces.
\newblock Diplomarbeit, Humboldt University, Berlin, 2002.

\bibitem[Wan89]{Wang:89}
McKenzy~Y. Wang.
\newblock Parallel spinors and parallel forms.
\newblock {\em Ann. Glob. Anal. and Geom.}, 7:59--68, 1989.

\bibitem[Wan93]{Wang:95}
McKenzy~Y. Wang.
\newblock On non-simply connected manifolds with non-trivial parallel spinors.
\newblock {\em Ann. Glob. Anal. and Geom.}, 13:31--42, 1993.

\end{thebibliography}

\vspace{0.8cm}
{\footnotesize Helga Baum, 
Institut f\"ur  Mathematik,
Humboldt-Universit\"at zu Berlin,
Sitz: Rudower Chaussee 25,
10099 Berlin,
Germany\\
email: baum@mathematik.hu-berlin.de\\}
{\footnotesize Ines Kath,
Max-Planck-Institut f\"ur Mathematik, Vivatsgasse 7, 53111 Bonn, Germany\\
email: kath@mpim-bonn.mpg.de\\}

\edo